\documentclass[a4paper,10pt]{amsart}

\usepackage{amssymb}
\usepackage[arrow]{xy}
\usepackage{pb-diagram,pb-xy}
\usepackage{graphicx,color,float}
\usepackage[a4paper,bookmarks]{hyperref}

\theoremstyle{plain}
\newtheorem{theorem}{Theorem}[section]
\newtheorem{proposition}[theorem]{Proposition}
\newtheorem{corollary}[theorem]{Corollary}
\newtheorem{lemma}[theorem]{Lemma}

\newtheorem*{folk}{Folklore Theorem}
\theoremstyle{definition}
\newtheorem{definition}[theorem]{Definition}

\newtheorem{example}[theorem]{Example}
\newtheorem{remark}[theorem]{Remark}
\newtheorem{notation}[theorem]{Notation}

\newtheorem*{thmb}{Theorem~\ref{theorem:homology-cobordism-of-homotopy-equivalent-manifolds}}
\newtheorem*{thmc}{Theorem~\ref{theorem:homology-invariance-der-series}}
\newtheorem*{thmd}{Corollary~\ref{corollary:homology-cobordism-spherical-space-form}}
\newtheorem*{thme}{Theorem~\ref{theorem:comparison}}
\newtheorem*{them}{Theorem}

\def\e{\epsilon}
\def\Z{\mathbb{Z}}
\def\Q{\mathbb{Q}}

\def\C{\mathbb{C}}
\def\L{\Lambda}

\def\GG{\mathcal{G}_\Gamma}
\def\O{\Omega}
\def\bO{\overline{\Omega}}

\def\HR{H\kern-0.1em R}
\def\HZ{H\mathbb{Z}}

\def\W{\mathcal{W}}

\def\MR{\mathcal{M}_R}

\def\N{\mathcal{N}}

\def\G{\mathcal{G}}
\def\K{\mathcal{K}}

\def\Ker{\operatorname{Ker}}
\def\Coker{\operatorname{Coker}}
\def\Im{\operatorname{Im}}
\def\Hom{\operatorname{Hom}}
\def\Tor{\operatorname{Tor}}

\def\sign{\operatorname{sign}}
\def\rank{\operatorname{rank}}

\def\inte{\operatorname{int}}

\def\ldim{\dim^{(2)}}
\def\lsign{\sign^{(2)}}
\def\L2{L^2}

\def\to{\mathchoice{\longrightarrow}{\rightarrow}{\rightarrow}{\rightarrow}}
\makeatletter
\newcommand{\shortxra}[2][]{\ext@arrow 0359\rightarrowfill@{#1}{#2}}
\def\longrightarrowfill@{\arrowfill@\relbar\relbar\longrightarrow}
\newcommand{\longxra}[2][]{\ext@arrow 0359\longrightarrowfill@{#1}{#2}}
\renewcommand{\xrightarrow}[2][]{\mathchoice{\longxra[#1]{#2}}%
  {\shortxra[#1]{#2}}{\shortxra[#1]{#2}}{\shortxra[#1]{#2}}}
\makeatother

\def\otimesover#1{\mathbin{\mathop{\otimes}_{#1}}}

\makeatletter
\def\Nopagebreak{\@nobreaktrue\nopagebreak}
\makeatother

\def\emptystr{}
\newcommand{\mkc}[2][]{\begin{color}{red}#2%
  \def\tempstr{#1}%
  \ifx\tempstr\emptystr \else\textsf{\SMALL\ \raise.7ex\hbox{[\tempstr]}}\fi
\end{color}}

\begin{document}

\title%
{$\L2$-signatures, homology localization, and amenable groups}

\author{Jae Choon Cha}

\address{Department of Mathematics and Pohang Mathematics Institute,
  Pohang University of Science and Technology, Pohang, Gyungbuk
  790--784, Republic of Korea}

\email{jccha@postech.ac.kr}

\author{Kent E. Orr}

\address{Department of Mathematics, Indiana University, Bloomington,
  Indiana 47405, USA}

\email{korr@indiana.edu}

\def\subjclassname{\textup{2000} Mathematics Subject Classification}
\expandafter\let\csname subjclassname@1991\endcsname=\subjclassname
\expandafter\let\csname subjclassname@2000\endcsname=\subjclassname


\begin{abstract}
  Aimed at geometric applications, we prove the homology cobordism
  invariance of the $\L2$-betti numbers and $\L2$-signature defects
  associated to the class of amenable groups lying in Strebel's
  class~$D(R)$, which includes some interesting infinite/finite
  non-torsion-free groups.  The proofs include the only prior known
  condition, that $\Gamma$ is a poly-torsion-free abelian group (or
  potentially a finite $p$-group.)
  We define a new commutator-type series which refines Harvey's
  torsion-free derived series of groups, using the localizations of
  groups and rings of Bousfield, Vogel, and Cohn.  The series, called
  the local derived series, has versions for homology with arbitrary
  coefficients, and satisfies functoriality and an injectivity
  theorem.  We combine these two new tools to give some applications
  to distinct homology cobordism types within the same simple homotopy
  type in higher dimensions, to concordance of knots in three
  manifolds, and to spherical space forms in dimension three.
\end{abstract}

\maketitle

\section{Introduction}

In their paper~\cite{Cochran-Orr-Teichner:1999-1}, Cochran, Orr, and
Teichner introduced $\L2$-signature defects (equivalently, von Neumann
$\rho$-invariants, or Cheeger-Gromov invariants) to study concordance
of knots, or more generally, homology cobordism classes of
$3$-manifolds.  They showed the invariance of $\L2$-signature defects
under integral homology cobordism over a poly-torsion-free abelian
group, that is, a group with a descending series admitting successive
$\Z$-torsion free quotients.  Since then, $\L2$-signatures have
appeared as a key ingredient of several interesting papers on homology
cobordism and concordance by authors including Cha, Cochran, Friedl,
Harvey, Heck, Horn, Kim, Leidy, Orr, and Teichner.

This paper substantially extends the above $\L2$-invariance under
homology cobordism to a much larger class of groups, and incorporates
homology with twisted coefficients as well.

Achieving this result requires a small shift of paradigm.  The
Cochran-Orr-Teichner results used the following property of a
poly-torsion-free abelian group~$\Gamma$: the group ring $\Z\Gamma$
embeds into a skew field which is a module over the Cohn localization
of~$\Z\Gamma$. This implies that a homology cobordism looks like a
product to the skew field.  For the groups $\Gamma$ we consider,
$\Z\Gamma$ may not embed in a skew field, requiring an entirely new
approach.  Our new approach, which uses neither Cohn localization nor
a skew field, subsumes poly-torsion-free abelian groups as a special
case.
We employ directly $\L2$-methods with coefficients in the von Neumann
algebra $\N\Gamma$ by using a result of L{\"u}ck on $\L2$-dimension
zero modules (see Theorem \ref{theorem:NG-is-$\L2$-dim-flat}).

We consider the class $D(R)$ of groups defined and studied by
Strebel~\cite{Strebel:1974-1}; for a commutative ring $R$ with unity,
a group $G$ lies in $D(R)$ if given any homomorphism of projective
$RG$-modules $\alpha\colon P\to Q$, then if $\alpha \otimes 1_R$
injective so is $\alpha$.  Our result below applies to the class of
\emph{amenable groups in $D(R)$}, a class which includes interesting
infinite and finite non-torsion-free groups, as well as, using the
case $R=\Z$, poly-torsion-free abelian groups (e.g., see
Lemma~\ref{lemma:amenable-and-D(R)}).  The following theorem gives a
special case of Theorem~\ref{theorem:homology-invariance-$\L2$-sign}.

\begin{them}  Suppose $G$ is an amenable group lying in Strebel's class
  $D(R)$.  If $W$ is an $R$-homology cobordism between two closed
  $n$-manifolds $M$ and $M'$ with restrictions to group homomorphisms
  $\phi$ and $\phi'$ of $\psi\colon \pi_1(W) \to G$ to $\pi_1(M)$ and
  $\pi_1(M')$, then the $\L2$-signature defects $\rho(M,\phi)$ and
  $\rho(M',\phi')$ are equal.
\end{them}

  To prove this, we use a new technique to control the $\L2$-dimension
  of homology with von Neumann group algebra $\N\Gamma$ coefficients,
  as shown in the following special case of
  Theorem~\ref{theorem:$\L2$-dim-local}.  This not only plays a key
  role in the proof of the above theorem, but we anticipate using this for future applications as well:
  
  \begin{them}
    Suppose $G$ is an amenable group in $D(R)$, and $C_*$ is a
    finitely generated free chain complex over~$\Z\Gamma$.  View $R$
    as a $\Z\Gamma$-module with trivial $\Gamma$-action.  If
    $H_*(C_*\otimes_{\Z\Gamma} R)=0$, then $H_*(C_*\otimes_{\Z\Gamma}
    \N\Gamma)$ has $\L2$-dimension zero.
  \end{them}

We apply our new theorem on $\L2$-signatures to study homology
cobordism classes of topological manifolds in dimension three and
higher, focusing on space forms, and three manifolds whose groups have
torsion elements.

We also employ an additional and essential new tool for these
results---a new commutator series, analogous to the Harvey derived
series of a group~\cite{Cochran-Harvey:2004-1}, but often much
smaller, allowing us to extract additional information from quotients.
We prove a type of Stallings injectivity theorem, similar in character
to those first developed by Cochran and Harvey
in~\cite{Cochran-Harvey:2004-1}, and proved similarly.  The
significance of this new series over the Harvey series is the use of
the Cohn localization of rings and modules, which yields a functorial
and often computable series, in place of the Ore localization used by
Harvey.  To define this series, we use group localization as well.  We
  call this new series, which can be defined for any coefficient $R$,
  the \emph{Vogel-Cohn $R$-local derived series}.  We define and
  investigate the analogous series using Bousfield localization as
  well.  (For more details, see
  Section~\ref{subsection:intro-new-series}
  and~\ref{subsection:local-derived-series}.)

As an application involving non-torsion-free groups, we give a
homology cobordism version of a theorem of Chang and
Weinberger~\cite{Chang-Weinberger:2003-1} on homeomorphism types of
manifolds with a given homotopy type.  In this paper we denote by
$\Z_{(p)}$ the classical localization of $\Z$ away from $p$, while
$\Z_p$ denotes $\Z/p\Z$.

\begin{thmb}
Suppose $M$ is a closed $(4k-1)$-manifold with $\pi=\pi_1(M)$, $k\ge
  2$.  Let $p$ be prime and $\pi^{(n)}$ be the $\Z_p$ or
  $\Z_{(p)}$-coefficient Vogel-Cohn local derived series of $\pi$.  If
  $\pi$ has a torsion element which remains nontrivial in
    $\pi/\pi^{(n)}$ for some $n$, then there exist infinitely many
  closed $(4k-1)$-manifolds $M_0=M$, $M_1$, $M_2,\ldots$ such that
  each $M_i$ is simple homotopy equivalent and tangentially equivalent
  to $M$ but $M_i$ and $M_j$ are not homology cobordant for any $i\ne
  j$.

\end{thmb}

In the proof, we make use of a nonvanishing property for certain
$\L2$-signatures associated to non-torsion-free groups due to Chang
and Weinberger~\cite{Chang-Weinberger:2003-1}, and apply our result to
capture the invariance of these $\L2$-signatures under homology
cobordism as well as homeomorphism.  The $\Z_p$ coefficient analogue
of
Theorem~\ref{theorem:homology-cobordism-of-homotopy-equivalent-manifolds}
holds as well (see
Section~\ref{subsection:distinct-homology-cobordism-types}.)

Additionally we apply these techniques to spherical 3-space forms and
prove the following.  We say that $M$ is \emph{homology equivalent to}
$N$ if there is a map $M\to N$ which induces isomorphisms on the
homology.

\begin{thmd}
  For any generalized quaternionic spherical 3-space form $M$, there
  are infinitely many closed 3-manifolds $M_0=M$, $M_1$, $M_2,\ldots$
  such that the $M_i$ are homology equivalent to $M$ and have
  identical Wall multisignatures (or equivalently Atiyah-Singer
  $G$-signatures) and Harvey $\L2$-signature invariants $\rho_n$
  \cite{Harvey:2006-1}, but no two of the $M_i$ are homology
  cobordant.
\end{thmd}

\subsection{Injectivity theorems}

The surprising injectivity theorem of Cochran and Harvey
\cite{Cochran-Harvey:2004-1} compares the commutator series quotients
of groups with similar homological properties, much as Stallings'
famous theorem did for the lower central series.  By an {\em
  injectivity theorem,} we mean a theorem of the type first proven by
Cochran and Harvey~\cite{Cochran-Harvey:2004-1}, and modeled after
Stallings' seminal work on the low dimensional homology of groups and
central series~\cite{Stallings:1965-1}.  Cochran and Harvey work with
a suitable commutator series of a group first suggested by Harvey, and
which we call here the Harvey series.  They show that every group
homomorphism $\pi \to G$ which induces an isomorphism $H_1(\pi, \Q)
\to H_1(G, \Q)$ and an epimorphism $H_2(\pi, \Q) \to H_2(G, \Q)$ also
induces a monomorphism modulo terms of the Harvey series.  The Harvey
series and associated Cochran-Harvey injectivity theorem underlie many
recent applications of the von Neumann $\rho$-invariants to
understanding concordance of knots and links, and homology cobordism
of manifolds.

\subsection{A new functorial commutator series and injectivity theorem}
\label{subsection:intro-new-series}

We consider the category $\GG$ of groups
$\pi$ over a fixed group~$\Gamma$, i.e., $\pi$ is endowed with a
homomorphism $\pi\to \Gamma$.  (In the special case $\Gamma=\{e\}$,
$\GG$ canonically identifies with the category of groups.)  In
Section~\ref{section:commutator-series}, for a given coefficient ring $R$
we define a new series
\[
\pi \supset \pi^{(0)} \supset \pi^{(1)} \supset \cdots \supset
\pi^{(n)} \supset \cdots
\]
of normal subgroups $\pi^{(n)}$ for each~$\pi \in \GG$ which we define
in terms of the Bousfield $R\Gamma$-homology localization of groups
and rings.  We call this $\{\pi^{(n)}\}$ the \emph{Bousfield
  $R\Gamma$-local derived series} to emphasize this is a functor on
the category $\GG$.  Similarly, we define the \emph{Vogel-Cohn
  $R\Gamma$-local derived series} using the Vogel localization of
groups and Cohn localization of rings.  Indeed the series can be
defined in a more general situation; for a precise description of
these local derived series, see
Definition~\ref{definition:local-derived-series}.

We prove that the $R\Gamma$-local derived series has the following
properties:

\begin{them}
  Let $\{\pi^{(n)}\}$ be the Bousfield (resp.\ Vogel-Cohn)
  $R\Gamma$-local derived series for $\pi$ in~$\GG$.
  \begin{enumerate}
  \item (Functoriality) For any morphism $\pi\to G$ in $\GG$, there
    are induced homomorphisms $\pi^{(n)} \to G^{(n)}$ and
    $\pi/\pi^{(n)} \to G/G^{(n)}$ for any~$n$.
  \item (Injectivity) If $\pi\to G$ is a group homomorphism which is
    2-connected on $H_*(-;R\Gamma)$ (resp.\ 2-connected on
    $H_*(-;R\Gamma)$ with $\pi$ finitely generated, $G$ finitely
    presented), then the induced map $\pi/\pi^{(n)} \to G/G^{(n)}$ is
    injective for any~$n$.
  \end{enumerate}
\end{them}

The above theorem combines portions of Lemma~\ref{lemma:naturality}
and Theorem~\ref{theorem:injectivity} from the body of the paper.

In Section~\ref{section:examples}, we give computational examples
which illustrate that our local derived series $\pi^{(n)}$ differs
from the Harvey series and other derived series discussed below.  In
particular, we illustrate that $\pi/\pi^{(n)}$ may have torsion and/or
infinite order elements, by comparison to the torsion free groups
$\pi/\pi_H^{(n)}$.

The above two advances, regarding local commutator series and
$\L2$-invariants, come together in the following theorem.

\begin{thmc}
  Let $R$ be either $\Z_p$, $\Z_{(p)}$, or $\Q$.  For a closed
  manifold $M$ over an amenable group $\Gamma$, view $\pi = \pi_1(M)$
  as a group over $\Gamma$ and denote the $R\Gamma$-coefficient
  Vogel-Cohn local derived series of $\pi$ over $\Gamma$ by
  $\{\pi^{(n)}\}$.  For the canonical map $\phi_n\colon \pi \to \pi /
  \pi^{(n)}$, the $\L2$-betti numbers $b^{(2)}_i(M,\phi_n)$ and the
  $\L2$-signature $\rho^{(2)}(M,\phi_n)$ are $R\Gamma$-homology
  cobordism invariants of~$M$ for any~$n<\infty$.  In particular, when
  $\Gamma$ is trivial, $b^{(2)}_i(M,\phi_n)$ and
  $\rho^{(2)}(M,\phi_n)$ are always $R$-homology cobordism invariants.
\end{thmc}

This applies to concordance of knots within a fixed homotopy class of
a three manifold.

\begin{corollary} (Compare to the Ph.D. thesis for Prudence
  Heck~\cite{Heck:2009-1}.)  Let $M$ be a three manifold with amenable
  fundamental group, let $g \in \Gamma = \pi_1(M)$.  Let $\pi =
  \pi_1(M-K)$ for some knot $K \subset M$ with $K$ in the homotopy
  class $g$, and let $\pi^{(n)}$ be the $n^{th} R\Gamma$-local derived
  series subgroup.  If $\pi^{(n)}/\pi^{(n+1)} \neq 0$, then there are
  infinitely many concordance classes of knots in the homotopy class
  of $g$ which for which no two are concordant to each other, and
  which are detected by the $\L2$-signature associated to $\pi
  \to \pi/\pi^{(n+1)}$.
\end{corollary}

Theorems concerning concordance of knots in general three manifolds
and using $L^2$-methods were first explored and proven in the
Ph.D. thesis of Prudence Heck~\cite{Heck:2009-1}.  The above argument
follows hers, with
Theorem~\ref{theorem:homology-invariance-der-series} applied at the
appropriate place in the argument, which the reader will find
apparent.  We refer the reader to this thesis~\cite{Heck:2009-1}.

\subsection{Comparison to other commutator series of groups}

Let $R$ be a subring of the rationals or a finite cyclic ring.  Recall
that for a group $\pi$, the (non-local) \emph{$R$-coefficient derived
  series,} or just \emph{$R$-derived series}, is defined recursively
as follows: $\pi^0 = \pi$, and given $\pi^n$, then
\[
\pi^{n+1} = \Ker{\left\{\pi^n \to \frac{\pi^n}{[\pi^n,
      \pi^n]} \otimesover{\Z} R\right\}}.
\]
Here, $[\pi, \pi]$ is the usual commutator subgroup of the group
$\pi$, the subgroup generated by commutators $[a, b] = aba^{-1}b^{-1}$
where $a$, $b \in \pi$.  When $R=\Z$, $\pi^n$ is the ordinary derived
series.  More generally, it is the fastest descending series such that
successive quotients are abelian and $R$-torsion free.  (This
definition of this series directly parallels the $R$-lower central
series of earlier authors, but was first applied
in~\cite{Harvey:2006-1}.)

To discuss alternative derived series of groups relevant in this work,
we first discuss group localization.

In~\cite{Bousfield:1975-1}, Bousfield constructed the $H_*(-; R)$
localization of a space.  Roughly, Bousfield associates to a space X,
a space $E(X)$ (up to homotopy type) and a map of spaces $X \to E(X)$
which is initial with the property that given an homology equivalence
$X\to Y$, there is a unique (homotopy class of a) map $Y \to E(X)$
making the following diagram homotopy commute.  (A more general
definition of localization is given in
Section~\ref{section:localization-groups-rings}.)

\[
\begin{diagram}
  \node{X} \arrow{se} \arrow[2]{e}
  \node[2]{Y}\arrow{sw,..}
  \\
  \node[2]{E(X)}
\end{diagram}
\]



Bousfield showed that the group $\pi_1(E(X))$ depends only on the group
$\pi_1(X)$.  Thus, one can define a functor on the category of groups
and group homomorphisms by $E(\pi) = \pi_1(E(B\pi))$, where $B\pi =
K(\pi, 1)$ is the classifying space of the discrete group~$\pi$.
Bousfield showed that $\pi \to E(\pi)$ localizes groups and group
homomorphisms with respect to $R$-homologically $2$-connected
homomorphisms; if $\pi\to G$ is $R$-homologically 2-connected, the
induced map $E(\pi)\to E(G)$ is an isomorphism.  Here $\pi \to G$ is
said to be \emph{$R$-homologically $2$-connected} if the induced
homomorphism $H_i(\pi; R) \to H_i(G; R)$ is an isomorphism for $i = 1$
and an epimorphism for $i = 2$.

Alternatively, in \cite{Vogel:1978-1} Vogel defined a localization
$\pi \to E(\pi)$ with respect to $R$-homologically 2-connected
homomorphisms $\pi\to G$ with $\pi$ finitely generated and $G$
finitely presented.  (See
Section~\ref{section:localization-groups-rings} for more details.)
Vogel developed this variation of Bousfield's work to study homology
types of compact manifolds, embeddings of compact manifolds, homology
cobordism, and the Cappell-Shaneson homology surgery groups.  (See
also~\cite{Cappell-Shaneson:1974-1}.)

The following series appears first in print in a paper of Cochran and
Harvey~\cite{Cochran-Harvey:2008-1} where they only consider Vogel
group localization.  In the definition below, $E(\pi)$ represents
either Vogel or Bousfield localization of groups, yielding distinct
functors.

\begin{definition}
\label{definition:R-local-derived-series}
  The \emph{universal $R$-local derived series of a group $\pi$},
  denoted using braces, as $\pi^{\{n\}}$, is recursively defined as
  follows: $\pi^{\{0\}} = \pi$.  Given the definition of
  $\pi^{\{n\}}$, define
  \[
  \pi^{\{n+1\}} = \Ker{\left\{\pi^{\{n\}} \to E(\pi)^{\{n\}} \to
      \frac{E(\pi)^{\{n\}}}{[E(\pi)^{\{n\}}, E(\pi)^{\{n\}}]}
      \otimesover{\Z} R\right\}}.
  \]
\end{definition}

We caution the reader that the universal $R$-local derived series does
not equal, in general, the $R$-derived series defined in a prior
paragraph.  Since $E^2=E$, one can see easily that $E(\pi)^{\{n\}} =
E(\pi)^n$ and $\pi^{\{n\}}$ is the preimage of $E(\pi)^n$ under
$\pi\to E(\pi)$.  In particular, the $R$-local and $R$-derived series
agree for local groups.

One easily proves the following analogue of Stallings'
theorem. Cochran and Harvey~\cite{Cochran-Harvey:2008-1} first
investigated this in depth, together with the above series, but the
theorem below appears to have been known already to experts in the
field.  In particular, Cochran and Harvey state the Vogel case in
terms of ``$R$-closures'' given in~\cite{Cha:2004-1}. (See
also~\cite{Levine:1989-1}.)  We label it a {\em folklore theorem}, but
wish to pay homage the aforementioned important contribution of
Cochran and Harvey.

\begin{folk}
  \leavevmode\Nopagebreak
  \begin{enumerate}
  \item (Bousfield case) Denote by $\pi^{\{n\}}$ the Bousfield universal
    $R$-local derived series.  Suppose $\alpha\colon \pi \to G$ is
    $R$-homologically 2-connected.  Then $\alpha$ induces an injection
    $\pi/\pi^{\{n\}} \to G/G^{\{n\}}$ for all~$n$.
  \item (Vogel case) Denote by $\pi^{\{n\}}$ the Vogel universal
    $R$-local derived series.  Suppose $\alpha\colon \pi \to G$ is
    $R$-homologically 2-connected with $\pi$ finitely generated and
    $G$ finitely presented.  Then $\alpha$ induces an injection
    $\pi/\pi^{\{n\}} \to G/G^{\{n\}}$ for all~$n$.
  \end{enumerate}
\end{folk}

\begin{proof}
  By construction $\pi/\pi^{\{n\}} \to E(\pi)/E(\pi)^{\{n\}}$ is
  injective for any $\pi$ and~$n$.  Since $E(\pi)=E(G)$, we have
  inclusions
  \[
  \pi/\pi^{\{n\}} \subset G/G^{\{n\}} \subset E(\pi)/E(\pi)^{\{n\}} =
  E(G)/E(G)^{\{n\}}. \qedhere
  \]
\end{proof}

The Harvey series~\cite{Cochran-Harvey:2004-1} may be thought of as a
rational coefficient approximation to the Vogel universal local
derived series, and the Cochran-Harvey injectivity theorem may be
viewed as a Harvey series version of this Folklore Theorem.  One issue
with the universal local derived series is that the series as
described is neither defined implicitly (that is it uses the
construction of $E(\pi)$ in its definition) nor is it readily
computable.  By contrast, the Harvey series consists of subgroups of
$\pi$ defined entirely from the structure of the group $\pi$ without
making use of the group localization $E(\pi)$, and one can compute
this in many useful circumstances.

On the other hand, a weakness of the Cochran-Harvey approach is that
their series and the resulting quotients are \emph{not} functorial
under homomorphisms of groups. That is, a group homomorphism $\pi \to
G$ does not induce $\pi^{(n)}_H \to G^{(n)}_H$ nor $\pi/\pi^{(n)}_H
\to G/G^{(n)}_H$ in general.

The $R\Gamma$-local series of this paper might be viewed as a
compromise theory --- one that lies closer to the ideal (i.e., the
universal local derived) series, that works for a wider range of
(twisted) coefficient systems, that is computable, that provides an
injectivity theorem, and importantly, that is functorial.

For $\Gamma = \{e\}$, the Vogel-Cohn local derived series refines
the Harvey series and lies closer to the universal local derived
series, as the following comparison theorem shows.

\begin{thme}
  Let $\Gamma=\{e\}$ and $\{\pi^{(n)}\}$ be the Vogel-Cohn $R\Gamma$-local
  derived series of a group~$\pi$.  Let $\{\pi^n\}$,
  $\{\pi^{\{n\}}\}$, $\{\pi^{(n)}_H\}$ be the $R$-derived series,
  Vogel universal $R$-local derived series, and the Harvey series,
  respectively.  Then, for any $\pi$ and $n$, we have
  \[
  \pi^n \subset \pi^{\{n\}} \subset \pi^{(n)} \subset  \pi^{(n)}_H.
  \]
\end{thme}

  In a future paper we will combine our results on homology with
  $\L2$-coefficients (Theorem~\ref{theorem:$\L2$-dim-local} and its
  applications), the universal $R\Gamma$-local derived series, and
  some computations using compact, orientable, three manifold groups
  whose group localization contains torsion, to obtain further
  examples of homology equivalent three manifolds which are not
  homology cobordant.

\subsection{Organization of the paper}

In the first half of this paper (Sections
\ref{section:localization-groups-rings}--\ref{section:examples}), we
define and investigate the $R\Gamma$-local derived series.  In
Section~\ref{section:localization-groups-rings}, we give necessary
preliminaries on localizations of groups and rings.  In
Section~\ref{section:commutator-series}, we define the $R\Gamma$-local derived
series and prove its functoriality and the injectivity theorem.  In
Section~\ref{section:comparison}, we compare the local derived series
with related commutator-type series, including the Harvey series.  In
Section~\ref{section:examples} we give some computational examples of
the local derived series which will be used in later applications.

In the latter half of this paper (Sections
\ref{section:$\L2$-dim-local-property}--\ref{section:applications}) we
study homological properties of $\L2$-theory and applications to
manifolds.  In Section~\ref{section:$\L2$-dim-local-property} we prove
our main theorem, the ``$\L2$-dimensional'' local property of the von
Neumann group ring of an amenable group, which play a crucial role in
our applications of the local derived series to manifolds.  {\em For
  those most interested in this result, one may read this section
  independently of the prior sections.}  In
Section~\ref{section:$\L2$-inv-homology-cobordism} we investigate
homology cobordism invariance of $\L2$-signature defects associated to
quotients of local derived series of fundamental groups.  In
Section~\ref{section:applications} we give some applications,
including
Theorems~\ref{theorem:homology-cobordism-of-homotopy-equivalent-manifolds}
Corollary~\ref{corollary:homology-cobordism-spherical-space-form}, and
some generalizations.

\subsection*{Acknowledgements}
This work was partly supported by a Korea Science and Engineering
Foundation grant funded by the Korean government (MOST) (R01--
2007--000--11687--0) and by a Korea Research Foundation grant funded
by the Korean Government (MOEHRD) (KRF--2007--412--J02302).  The second author was supported by NSF grant DMS-0707078.  The first
author thanks the Department of Mathematics, Indiana University, for
supporting his visit in 2008 summer.

\section{Localizations of groups and rings}
\label{section:localization-groups-rings}

In this section we discuss basic definitions and properties of
homology localizations of groups and rings, which are essentially due
to Bousfield, Vogel, and Cohn.  Theorems and propositions in this
section are not new and have been widely known (at least to experts in
related fields).  We will focus on those facts that we will need in
later sections.

We begin with the general definition of a localization functor in a
category.

\begin{definition}
  \label{definition:localization}
  Suppose $\O$ is a collection of morphisms in a
  category~$\mathcal{C}$.
  \begin{enumerate}
  \item An object $A$ in $\mathcal{C}$ is \emph{local with respect to
      $\O$} if for any $\pi \to G$ in $\Omega$ and for any $\pi \to
    A$, there exists a unique morphism $G \to A$ making the following
    diagram commute:
    \[
    \begin{diagram}
      \node{\pi}\arrow{e}\arrow{s}
      \node{G} \arrow{sw,..} \\
      \node{A}
    \end{diagram}
    \]

  \item A \emph{localization functor} with respect to $\O$ is a
    functor 
    $$
    E\colon \mathcal{C} \to \{\mbox{local objects in $\mathcal{C}$}\}
    $$ 
    endowed with a natural transformation $\{p_G\colon
    G \to E(G)\}_{G\in\mathcal{C}}$ such that for any morphism $G\to
    A$ into a local object $A$, there is a unique morphism $E(G) \to
    A$ making the following diagram commute:
    \[
    \begin{diagram}
      \node{G}\arrow{e,t}{p_G}\arrow{s}
      \node{E(G)} \arrow{sw,..} \\
      \node{A}
    \end{diagram}
    \]
  \end{enumerate}
\end{definition}

The following statements are consequences of the universal properties
stated above.  We omit proofs.
\begin{proposition}
  \leavevmode\nolinebreak
  \begin{enumerate}
  \item A localization functor is unique (up to a natural
    equivalence).
  \item A localization functor $E$ is an idempotent, i.e., $E^2 = E$.
  \item If E is a localization functor with respect to $\Omega$, then
    every $\pi\to G$ in $\Omega$ induces an equivalence $E(\pi)\to E(G)$.
  \end{enumerate}
\end{proposition}

\subsection{Localizations of groups}
We will give a very brief introduction to a homology localization
theory of groups due to Bousfield and Vogel.  Let $\GG$ be the
category of groups over a fixed group $\Gamma$.  Precisely, the
objects of $\GG$ are homomorphisms of a group $\pi$ into~$\Gamma$.  A
morphism from $\pi \to \Gamma$ to $G \to \Gamma$ is a homomorphism
$\pi\to G$ making the following diagram commute:
\[
\begin{diagram}
  \node{\pi} \arrow{se}\arrow[2]{e}
  \node[2]{G} \arrow{sw}
  \\
  \node[2]{\Gamma}
\end{diagram}
\]
We abuse notation and denote an object $\pi\to \Gamma$ of
$\GG$ by~$\pi$.

As a special case, when $\Gamma$ is a trivial group, $\GG$ is
canonically equivalent to the category of groups.

Let $R$ be either a finite cyclic ring $\Z_d$ or a subring of $\Q$.
For an object $\pi$ in $\GG$, the homology $H_*(\pi;R\Gamma)$ with
local coefficients is defined, where $R\Gamma$ denotes the group ring
of $\Gamma$ over~$R$.  We will consider localizations with respect to
$R\Gamma$-coefficient homology.  Specifically, let $\HR$ be the collection of morphisms
$\alpha\colon \pi \to G$ in $\GG$ which are 2-connected on
$H_*(-;R\Gamma)$, that is, $\alpha$ induces an isomorphism on
$H_1(-;R\Gamma)$ and an epimorphism on $H_2(-;R\Gamma)$.  A
localization functor on $\GG$ with respect to $\HR$ is called the
\emph{Bousfield $\HR$-localization}.

Regarding applications to geometric topology (especially to the study
of compact manifolds), one is naturally led to consider the
subcollection of morphisms on finitely presented groups.  Motivated
by this, we also consider the subcollection $\Omega_0$ of morphisms
$\alpha\colon \pi\to G$ in $\HR$ with $\pi$ and $G$ finitely
presented.  A localization functor on $\GG$ with respect to $\O_0$ is
called the \emph{Vogel localization}.

The following results are essentially due to
Bousfield~\cite{Bousfield:1974-1,Bousfield:1975-1} and
Vogel~\cite{Vogel:1978-1}.  (See also
Farjoun-Orr-Shelah~\cite{Farjoun-Orr-Shelah:1989},
Levine~\cite{Levine:1989-1}, and Cha~\cite{Cha:2004-1} for
equation-based approaches.  Levine appears to have first observed
these combinatorial approaches within the context of Vogel localization, with Farjoun and Shelah doing so later
but independently in the context of Bousfield's work.)

\begin{theorem}[Bousfield~\cite{Bousfield:1975-1,Bousfield:1974-1},
  Vogel~\cite{Vogel:1978-1}]
  \label{theorem:existence-of-group-localization}
  For any $R$ and $\Gamma$, there exist a Bousfield
  $\HR$-localization functor and a Vogel localization functor on~$\GG$.
\end{theorem}

In the case of the Vogel localization $E$, the map $\pi \to E(\pi)$ is
not always contained in $\O_0$, even when $\pi$ is finitely presented.
In this case, one can enlarge $\O_0$ to a collection, $\bO_0$ which
contains $\pi \to E(\pi)$, and for which the resulting localization
functors agree on classes in $\O_0$.  Let $\bO_0$ be the collection of
morphisms $\pi\to G$ such that for any given commutative diagram
\[
\begin{diagram}
  \node{A} \arrow{e} \arrow{s}
  \node{B} \arrow{s}
  \\
  \node{\pi} \arrow{e}
  \node{G}
\end{diagram}
\]
of morphisms in $\GG$ with $A$ and $B$ are finitely presented, there
is $A_0 \to B_0$ in $\O_0$ which fits into the following commutative
diagram:
\[
\begin{diagram}\dgARROWLENGTH=0.6\dgARROWLENGTH
  \node{A} \arrow[3]{e} \arrow[2]{s} \arrow{se,..}
  \node[3]{B} \arrow[2]{s} \arrow{sw,..}
  \\
  \node[2]{A_0} \arrow{e} \arrow{sw,..}
  \node{B_0} \arrow{se,..}
  \\
  \node{\pi} \arrow[3]{e}
  \node[3]{G}
\end{diagram}
\]

\begin{remark}
  \label{remark:bar-omega-is-in-HR}
  It is easily seen that $\bO_0 \subset \HR$, since any homology class
  of a free chain complex is supported by a finitely generated
  subcomplex.
\end{remark}

\begin{theorem}[Bousfield~\cite{Bousfield:1975-1,Bousfield:1974-1},
  Vogel \cite{Vogel:1978-1}]
  \label{theorem:properties-of-group-localization}
  \leavevmode\nolinebreak
  \begin{enumerate}
  \item For the Bousfield $\HR$-localization $E$, $\pi\to E(\pi)$
    is in~$\HR$ for all $\pi\in \GG$.
  \item The Vogel localization $E$ is a localization with respect to
    $\bO_0$, and in this case, $\pi \to E(\pi)$ is in $\bO_0$ for all
    $\pi\in\GG$.
  \end{enumerate}
\end{theorem}

We sketch the proofs of
Theorems~\ref{theorem:existence-of-group-localization} and
\ref{theorem:properties-of-group-localization} in an appendix for the
convenience of our readers.

\begin{remark}
  \label{remark:normal-generation-condition}
  In known applications of the localization theory to the study of
  manifold embeddings, one is naturally led to consider localizations
  with respect to certain morphisms $\alpha\colon\pi \to G$ in $\O$
  with the property that $\Ker\{G\to\Gamma\}$ is normally generated by
  the image of $\Ker\{\pi\to\Gamma\}$.  (See, for instance,
  \cite{LeDimet:1988-1}.)  It is known that
  Theorems~\ref{theorem:existence-of-group-localization} and
  \ref{theorem:properties-of-group-localization} also hold for the
  analogues of $\HR$, $\O_c$ and $\bO_c$ with this normal generation
  condition.  It can also be seen that all results in this paper hold
  under the normal generation condition.
\end{remark}

 \begin{remark}
  \label{remark:localization-with-cardinality-condition}
  Following Vogel, one may also consider only homomorphisms $\pi \to
  G$ in $\HR$ such that the $\pi$ and $G$ have a number of generators
  and relations each bounded below by a fixed cardinal, $c$.  For
  instance, the Vogel localization of groups is localization with
  respect to $\HR$ homomorphisms and the cardinal $c=\aleph_0$.
  However, one can show that for any cardinal $c > \aleph_0$, the
  resulting localization equals the usual Bousfield $\HR$
  localization.  Analogous observations apply equally to the (module
  and) ring localizations which will de discussed below.
\end{remark}

\subsection{Localizations of rings}
\label{subsection:ring-localization}

We begin by recalling Cohn's classical ring localization.  Let $R$ and
$S$ be rings with unity, together with a ring homomorphism $\epsilon :
R \to S$.

The \emph{Cohn localization}, $\Lambda$, of the homomorphism $\e\colon
R \to S$ is a ring endowed with a ring homomorphism $R\to \Lambda$
satisfying the following conditions:
\begin{enumerate}
\item Given any square matrix $A$ over $R$, then $A \otimes_R S$ in
  invertible if and only of $A \otimes_{R} \Lambda $ is invertible,
  and
\item 
  $R\to \Lambda$ is initial among homomorphisms of $R$ satisfying~(1),
  that is, if $R\to \Lambda'$ is another ring homomorphism
  satisfying~(1), then there is a unique ring homomorphism $\Lambda
  \to \Lambda'$ making the following diagram commute:
  \[
  \begin{diagram}
    \node{R}\arrow[2]{e} \arrow{se}
    \node[2]{\Lambda} \arrow{sw,..}
    \\
    \node[2]{\Lambda'}
  \end{diagram}
  \]
\end{enumerate}

In~\cite{Cohn:1971-1}, Cohn constructed this localization and showed
uniqueness and functoriality up to isomorphism.  We give an
alternative description of Cohn localization here as a localization of
modules with a naturally imposed ring structure.  (Of particular
interest in this paper is the case of $R\pi \to R\Gamma$ induced by a
homomorphism $\pi\to\Gamma$ in $\GG$.)

Let $\W$ be the collection of homomorphisms $\alpha\colon F \to F'$
between (right) free $R$-modules $F$ and $F'$ with the same rank such
that
\[
\alpha\otimes 1_{S} \colon F\otimes_{R}S \to
F'\otimes_{R}S
\]
is an isomorphism.  (Note that the rank of $F$ and $F'$ may be an
arbitrary cardinal.)  Recall that, with respect to $\W$, a local
  $R$-module is defined to be an $R$-module such that for any
  $\alpha\colon F\to F'$ in $\W$,
\[
\alpha^*\colon \Hom{(F', M)} \to \Hom{(F, M)}
\]
is an isomorphism.  A localization functor with respect to $\W$
\[
E\colon \{\text{$R$-modules}\} \to \mathcal \{\text{local
  $R$-modules}\}
\] 
in the sense of Definition~\ref{definition:localization} is called the
\emph{Bousfield $\HZ$-localization functor}.  (See
Remark~\ref{remark:original-Bousfield-HZ}.)

We also consider module localization with an additional finite rank
condition: let $\W_0$ be the collection of $\alpha\colon F\to F'$ in
$\W$ such that $\rank F = \rank F'$ is finite.  We call the resulting
localization functor of $R$-modules the \emph{Cohn localization}.  The
following existence results are essentially due to Bousfield and Cohn:

\begin{theorem}[Bousfield~\cite{Bousfield:1975-1},
  Cohn~\cite{Cohn:1971-1}]
  For any $R \to S$, there exist the Bousfield
  $\HZ$-localization and Cohn localization of modules.
\end{theorem}

The Bousfield $\HZ$ and Cohn localization of modules have the
  following properties:

\begin{enumerate}
\item[($1$)] If $E$ represents either the Bousfield $\HZ$ and Cohn
  localization functor of $R$-modules, then the $R$-module $E(R)$ admits a
  natural ring structure such that $R \to E(R)$ is a ring
  homomorphism.
\item[($2$)] The Bousfield $\HZ$-localization $E(R)$ has the following
  universal property as a ring: let $\mathcal{C}^\W$ be the category
  of rings $\Lambda$ endowed with a ring homomorphism $R\to\Lambda$
  such that for any $\alpha\colon F\to F'$ in $\W$,
  \[
  \alpha\otimes 1_\Lambda\colon F\otimes_R \Lambda \to F'\otimes_R\Lambda
  \]
  is an isomorphism.  Then $E(R)$ is the initial object of
  $\mathcal{C}^\W$.
\item[($3)$] The Cohn module localization of $R$ with the ring structure from (1) has the universal
  property analogous to ($2$) where $\W_0$ plays the role of~$\W$.
\end{enumerate}

We give proofs in an appendix.  Property ($3$) says that the
Cohn localization of the $R$-module $R$ agrees with the aforementioned
Cohn localization of the ring $R$ (endowed with $\epsilon\colon R\to
S$).  In parallel to this, we call $E(R)$ in ($2$) the
\emph{Bousfield $\HZ$-localization of the ring $R$} endowed with
$\epsilon\colon R \to S$.  We note that this viewpoint enables us to
regard both ring localizations as special cases of
Definition~\ref{definition:localization} and to use standard
properties of a localization functor.

\begin{remark}\label{remark:original-Bousfield-HZ}
  In~\cite{Bousfield:1975-1}, Bousfield originally defined
  $\HZ$-localization of $R\pi$-modules for $R\pi \to R\Gamma$ induced
  by a group homomorphism $\pi\to\Gamma$.  He localizes with respect
  to the class of $R\pi$-module homomorphisms $\alpha\colon A\to B$
  such that $\Tor_i^{R\pi}(R\Gamma,A) \to \Tor_i^{R\pi}(R\Gamma,B)$ is
  an isomorphism for $i=0$ and a surjection for $i=1$.  He dealt
  mainly with the case of $\Gamma=\{e\}$ and $R=\Z$, so that the terms
  $\Tor_i^{R\pi}(R\Gamma,-)$ equal $H_i(\pi;-)$, but his arguments
  extend to the case of any $\Gamma$ and~$R$.  This module
  localization agrees with the module localization with respect to
  $\W$ described above.  We prove this in an appendix.
  \end{remark}

We note that by applying the universal property of the Bousfield $\HZ$
(or Cohn) localization $\Lambda$ for $\epsilon\colon R\to S$, there is
a unique homomorphism $\Lambda\to S$ such that the following diagram
commutes:
\[
\begin{diagram}
  \node{R} \arrow[2]{e} \arrow{se} \node[2]{\Lambda} \arrow{sw,..}
  \\
  \node[2]{S}
\end{diagram}
\]
The Bousfield $\HZ$ (resp.\ Cohn) localization is the initial object
of the subcategory of objects of $\mathcal{C}^{\W}$ (resp.\
$\mathcal{C}^{\W_0}$) admitting such a diagram.

\begin{notation}
  Given $\pi \to \Gamma$ an element of the category $\GG$, we denote
  the Cohn and Bousfield $\HZ$-localization of $R\pi$ (endowed with
  $R\pi \to R\Gamma$) by $L(\pi)=\Lambda$ and observe that this is
  functorial on morphisms in $\GG$.
\end{notation}

\subsection{Ring localizations and homology}
The following property of the ring localizations will play a key
role later.  It is essentially due to Vogel.

\begin{theorem}[Vogel~\cite{Vogel:1982-1}]
  \label{theorem:2-conn-prop-of-ring-localization}
  Suppose $\pi\to G$ is a morphism in~$\GG$.
  \begin{enumerate}
  \item The Bousfield $\HZ$-localization $\Lambda$ of $RG$ satisfies
    the following: every $\pi \to G$ in~$\HR$ is 2-connected on
    $H_*(-;\Lambda)$.
  \item The Cohn localization $\Lambda$ of $RG$ satisfies the
    following: every $\pi \to G$ in $\bO_{0}$ is 2-connected on
    $H_*(-;\Lambda)$.
  \end{enumerate}
\end{theorem}

\begin{proof}
  We first prove (2).  Suppose $\pi\to G$ is in~$\bO_{0}$.  Let
  $C_*(G,\pi;\Lambda)$ be the $\Lambda$-coefficient CW chain complex
  of the mapping cylinder $M$ of $K(\pi,1) \to K(G,1)$.  All
  CW-complexes are over $K(\Gamma,1)$.

  For a cycle $z$ in $C_i(G,\pi;\Lambda)$ ($i=1,2$), there is a finite
  subpair of $(M,K(\pi,1))$ in which $z$ is supported.  Therefore, we
  obtain the following commutative diagram in $\GG$
  \[
  \begin{diagram}
    \node{A} \arrow{s} \arrow{e}
    \node{B} \arrow{s}
    \\
    \node{\pi} \arrow{e}
    \node{G}
  \end{diagram}
  \]
  such that $A$ and $B$ are the finitley presented fundamental groups
  of this subpair and $z$ is the image of some cycle in
  $C_*(B,A;\Lambda)$.  Since $\pi\to G$ is in $\bO_0$, we obtain
  \[
  \begin{diagram}\dgARROWLENGTH=0.6\dgARROWLENGTH
    \node{A} \arrow[3]{e} \arrow[2]{s} \arrow{se,..}
    \node[3]{B} \arrow[2]{s} \arrow{sw,..}
    \\
    \node[2]{A_0} \arrow{e} \arrow{sw,..}
    \node{B_0} \arrow{se,..}
    \\
    \node{\pi} \arrow[3]{e}
    \node[3]{G}
  \end{diagram}
  \]
  where $A_0 \to B_0$ is a morphism in $\Omega_0$.

  Since $z$ is the image of a cycle in $C_*(B_0,A_0;\Lambda)$, in
  order to conclude $z$ bounds it suffices to prove the following
  claim: $H_i(B_0, A_0;\Lambda)=0$ for $i=1,2$.  Since
  $C_i(B_0,A_0;RG)$ can be assumed to be $RG$-free module of finite
  rank for $i\le 2$, the claim is proved by a standard partial chain
  contraction argument, which is originally due to Vogel.  (For
  example, refer to \cite{Vogel:1982-1}, \cite{Levine:1994-1},
  \cite{Cochran-Orr-Teichner:1999-1}.)  This completes the proof of
  (2).

  The proof of (1) is similar but easier since we do not have to
  consider $A\to B$ and $A_0\to B_0$; we can directly apply the above
  chain contraction argument to $C_*(G, \pi; R\pi)$ to show that
  $H_i(G, \pi; \Lambda)=0$ for $i=1,2$.
%
\end{proof}

\section{Local derived series and injectivity}
\label{section:commutator-series}

Fix $R$ and let $\GG$ be as before.  In this section we consider a
collection $\O$ of morphisms in $\GG$, a localization functor $E$ on
$\GG$, and a localization functor~$L$ for group rings; suppose $(\O,E,L)$
is either one of the following two cases:
\[
\text{(i) }\left\{
  \begin{aligned}
    \O &= \HR \\
    E &= \text{Bousfield $\HR$-localization} \\
    L &= \text{Bousfield $\HZ$-localization}
  \end{aligned}
\right\}
\text{ or }
\text{(ii) }\left\{
  \begin{aligned}
    \O &= \bO_0 \\
    E &= \text{Vogel localization} \\
    L &= \text{Cohn localization}
  \end{aligned}
\right\}
\]  

We will define a commutator-type series for each $\pi \in \GG$ which
is determined by $R$ and $(\O, E, L)$ and prove that it admits an
injectivity theorem.  Indeed all the properties of $(\O,E,L)$ we need
are the following, which are easily verified for the above cases (i)
and (ii), by our discussion in the previous section (in particular see
Theorems~\ref{theorem:properties-of-group-localization}
and~\ref{theorem:2-conn-prop-of-ring-localization}):

\begin{enumerate}
\item $E$ is a localization functor on $\GG$ with respect to $\O$.
\item $\pi\to E(\pi)$ is in $\Omega$ for any $\pi\in \GG$.
\item Every $\pi\to G$ in $\O$ is 2-connected on $H_*(-;\Lambda)$,
  where $\Lambda=L(G)$.
\end{enumerate}

\begin{definition} We call $(\O,E,L)$ a \emph{homology localization
    triple} if the above (1), (2), and (3) hold.  We call the homology
  localization triples (i) and (ii) described above the
  \emph{Bousfield localization triple} and \emph{Vogel-Cohn
    localization triple}, respectively.
\end{definition}

Throughout this section, $(\O,E,L)$ is always assumed to be a homology
localization triple, and $\widehat\pi$ denotes $E(\pi)$ for $\pi\in
\GG$.

\subsection{$R\Gamma$-local derived series}
\label{subsection:local-derived-series}

\begin{definition}
\label{definition:local-derived-series}
  For each $\pi$ in $\GG$, $(\O,E,L)$ a homology localization triple,
  and for each ordinal $n$, the $R\Gamma$-\emph{local derived series}
  $\{\pi^{(n)}\}$ is defined transfinite-inductively as follows:
  \[
  \pi^{(0)} = \Ker\{\pi \to \Gamma\},
  \]
  and assuming $\pi^{(n)}$ has been defined, $\pi^{(n+1)}$ is defined by
  \[
  \pi^{(n+1)} = \Ker \bigg\{ \pi^{(n)} \to \frac{\pi^{(n)}}{[\pi^{(n)},
    \pi^{(n)}]} = H_1\Big(\pi;\Z\Big[\frac{\pi}{\pi^{(n)}}\Big]\Big)
  \to H_1(\pi;\Lambda_\pi) \bigg\}
  \]
  where $\Lambda_\pi$ designates the ring localization
  $L(\widehat\pi/\widehat\pi^{(n)})$, and recall that $\widehat \pi =
  E(\pi)$.  For a limit ordinal $n$, define $\pi^{(n)}=\bigcap_{k<n}
  \pi^{(k)}$.
\end{definition}

By an induction one easily verifies that each $\pi^{(n)}$ is a normal
subgroup of~$\pi$ forming a normal series
\[
\pi \supset \pi^{(0)} \supset \pi^{(1)} \supset \cdots \supset
\pi^{(n)} \supset \cdots.
\]

\begin{lemma}
  \label{lemma:naturality}
  For any $n$, the associations $\pi \to \pi^{(n)}$, $\pi \to
  \pi/\pi^{(n)}$ are functors $\GG \to \GG$.  In particular, any
  morphism $\pi \to G$ in $\GG$ induces group homomorphisms $\pi^{(n)}
  \to G^{(n)}$ and $\pi/\pi^{(n)} \to G/G^{(n)}$.
\end{lemma}

\begin{proof}
  We use a transfinite induction on $n$ to show that a morphism $\pi
  \to G$ gives rise to a group homomorphism on the normal subgroups
  $(-)^{(n)}$ and so on their quotients.  For $n=0$, the conclusion is
  obvious.  Assume that the conclusion holds for~$n$.  Let
  $\Lambda_\pi=L(\widehat\pi/\widehat\pi^{(n)})$ and
  $\Lambda_G=L(\widehat G/\widehat G^{(n)})$ as before.  Then since
  $L$ is functorial, there is an induced map $\Lambda_\pi \to
  \Lambda_G$.  So we have a commutative diagram
  \[
  \begin{diagram}
    \node{\pi^{(n)}} \arrow{s}\arrow[2]{e}
    \node[2]{H_1\big(\pi ; \Lambda_\pi\big)}
    \arrow{s}
    \\
    \node{G^{(n)}} \arrow[2]{e}
    \node[2]{H_1\big(G;\Lambda_G\big)}
  \end{diagram}
  \]
  where the horizontal maps are the compositions appearing in the
  definition of $\pi^{(n+1)}$ and $G^{(n+1)}$.  From this it follows
  that $\pi^{(n+1)}$ is sent into~$G^{(n+1)}$.

  For a limit ordinal $n$, assuming the conclusion holds for $k<n$,
  the conclusion for $n$ follows immediately by taking the limit.
\end{proof}

\begin{remark}
  The above naturality is essentially due to the functoriality of $L$
  and $E$ used in the above definition. By contrast, the
  Harvey series $\{\pi_H^{(n)}\}$
  defined in~\cite{Cochran-Harvey:2004-1} is not functorial on the
  category of groups since the Ore localization of group rings
  is not functorial.
\end{remark}

The following Lemma is proven easily from the
definition of the local derived series.

\begin{lemma}
  Suppose $(\O,E,L)$ and $(\O',E',L')$ are homology localization
  triples (possibly for distinct $R$) and there are natural
  transformations $E\to E'$ and $L\to L'$.  Then the local derived
  series determined by $(\O,E,L)$ is contained in the local derived
  series determined by $(\O',E',L')$.
\end{lemma}

\begin{definition}
  \label{definition:Vogel-Cohn-Bousfield-local-derived-series}
  We have the following special cases of
  Definition~\ref{definition:local-derived-series} using the Bousfield
  and Vogel-Cohn homology localization triples, cases (i) (resp.\
  (ii)) described in the beginning of
  section~\ref{section:commutator-series}.  We denote this series
  $\{\pi^{(n)}\}$ and call it the \emph{Bousfield (resp.\ Vogel-Cohn)
    $R\Gamma$-local derived series of $\pi$ over $\Gamma$}.
\end{definition}

\begin{proposition}
  \leavevmode\nolinebreak
  \begin{enumerate}
  \item If $R\subset R'$, then the $R$-coefficient Bousfield local
    derived series is contained in the $R'$-coefficient Bousfield
    local derived series.  Similarly for the Vogel-Cohn case.
  \item For any $R\subset \Q$, the Vogel-Cohn local derived series is
    contained in the Bousfield local derived series.
  \end{enumerate}
\end{proposition}

The existence of the required natural transformations follows from the
universal property of a localization and the facts $\O_0 \subset \HR$
and $\W_0 \subset \W$.  We omit details.

\subsection{An injectivity theorem}

\begin{theorem}
  \label{theorem:injectivity}
  Suppose $(\O,E,L)$ is a homology localization triple and let
  $\{\pi^{(n)}\}$ be the associated local derived series.  Then for
  any $\pi \to G$ in $\Omega$ and for any ordinal $n$, the induced
  homomorphism $\pi/\pi^{(n)} \to G/G^{(n)}$ is injective.  In
  particular, $\pi/\pi^{(n)} \to \widehat\pi / \widehat\pi^{(n)}$ is
  injective.
\end{theorem}

\begin{proof}
  We use a transfinite induction on~$n$ to prove the first statement.
  Once we prove this, the second statement follows immediately since
  $\pi\to \widehat\pi$ is in~$\O$.

  Since
  \[
  \begin{diagram}
    \node{\pi} \arrow[2]{e} \arrow{se}
    \node[2]{G} \arrow{sw}
    \\
    \node[2]{\Gamma}
  \end{diagram}
  \]
  commutes, $\pi^{(0)}= \Ker\{\pi\to \Gamma\}$ is exactly the inverse
  image of $G^{(0)}=\Ker\{G \to \Gamma\}$.  Therefore the conclusion
  holds for $n=0$.

  Suppose it holds for~$n$.  Denote $L(\widehat\pi/\widehat\pi^{(n)})$
  by~$\Lambda_\pi$ as before.  To show the conclusion for $n+1$,
  consider the commutative diagram
  \[
  \begin{diagram}\dgARROWLENGTH=1.5em
    \node{1} \arrow{e}
    \node{\pi^{(n)}/\pi^{(n+1)}} \arrow{e}\arrow{s}
    \node{\pi/\pi^{(n+1)}} \arrow{e}\arrow{s}
    \node{\pi/\pi^{(n)}} \arrow{e}\arrow{s}
    \node{1}
    \\
    \node{1} \arrow{e}
    \node{G^{(n)}/G^{(n+1)}} \arrow{e}
    \node{G/G^{(n+1)}} \arrow{e}
    \node{G/G^{(n)}} \arrow{e}
    \node{1}
  \end{diagram}
  \]
  with exact rows.  By the five lemma and the induction hypothesis, it
  suffices to show the injectivity of the homomorphism
  $\pi^{(n)}/\pi^{(n+1)} \to \G^{(n)}/G^{(n+1)}$.

  Since $\pi\to G$ is in $\Omega$, we have $\widehat\pi \cong \widehat
  G$ and $\Lambda_\pi \cong \Lambda_G$.  Therefore, we have the
  following commutative diagram, where $\pi^{(n+1)}$ and $G^{(n+1)}$
  are the kernels of the horizontal arrows:
  \[
  \begin{diagram}
    \node{\pi^{(n)}} \arrow{s}\arrow[2]{e} \node[2]{H_1(\pi;
      \Lambda_\pi)} \arrow{s}
    \\
    \node{G^{(n)}} \arrow[2]{e} \node[2]{H_1(G;\Lambda_\pi)}
  \end{diagram}
  \]
  By the defining properties of a homology localization triple,
  $\pi\to G$ is 2-connected on $H_*(-;\Lambda)$ where
  $\Lambda=L(\pi)$.  Since there is a natural homomorphism $\Lambda \to
  \Lambda_\pi$ by the naturality of the ring localization $L$, $\pi\to
  G$ is also 2-connected on $H_*(-;\Lambda_\pi)$.  Consequently the
  rightmost vertical arrow in the above diagram is an isomorphism.
  From this it follows that the inverse image of $G^{(n+1)}\subset
  G^{(n)}$ under $\pi^{(n)} \to G^{(n)}$ is exactly $\pi^{(n+1)}$.
  This proves the conclusion for~$n+1$.  The limit ordinal case
  follows by taking limits.
\end{proof}

\begin{remark}
  \leavevmode\Nopagebreak
  \begin{enumerate}
  \item Consider the case of the rational Bousfield or Vogel-Cohn
    derived series over $\Gamma=\{e\}$.  Then as a part of the
    hypotheses of Theorem~\ref{theorem:injectivity}, we have an
    $H_1$-condition that $\pi\to G$ induces an isomorphism on
    $H_1(-;\Q)$.  This can be weakened to only require that the
    induced map on $H_1(-;\Q)$ be injective.  For, one can think of an
    appropriate map from the free product $\pi * F$ with a free group
    $F$ into $G$ which is 2-connected, so that $(\pi*F)/(\pi*F)^{(n)}
    \to G/G^{(n)}$ is injective by Theorem~\ref{theorem:injectivity}.
    Since $\pi\to \pi*F$ has a left inverse, so does $\pi/\pi^{(n)}\to
    (\pi*F)/(\pi*F)^{(n)}$ by functoriality.  It follows that
    $\pi/\pi^{(n)} \to G/G^{(n)}$ is injective.

  \item For the Vogel-Cohn derived series over $\Gamma=\{e\}$ (for any
    $R$), one can prove the following Dwyer-type injectivity theorem,
    which is similar to a result in~\cite{Cochran-Harvey:2006-01}: if
    a homomorphism $\pi \to G$ between finitely presented groups
    induces an isomorphism $\widehat\pi \cong \widehat G$ and an
    epimorphism
    \[
    H_2(\pi;R) \to H_2(G)/\Im\big\{ H_2(G^{(n)};R) \to H_2(G;R)\big\},
    \]
    then $\pi/\pi^{(n+1)} \to G/G^{(n+1)}$ is injective.  The proof is
    omitted.
  \end{enumerate}
\end{remark}

\section{Comparison with Harvey series}
\label{section:comparison}

\def\univ{\text{U}}
\def\VC{\text{VC}}
\def\Harvey{H}

In this section, for $\Gamma = \{e\}$, we compare the $R\Gamma$-local
derived series with the Harvey series~\cite{Cochran-Harvey:2004-1} and
other related series.  Because we need to distinguish several series,
we use the following notation (temporarily): Fix a coefficient ring
$R$, and for a group $\pi$, let $\{\pi^{n}\}$, $\{\pi^{\{n\}}\}$, and
$\{\pi^{(n)}\}$ be the $R$-coefficient derived series (defined
immediately above Definition~\ref{section:commutator-series}),
universal local derived series (see
Defintion~\ref{section:commutator-series}), and Vogel-Cohn $R$-local
derived series (see
Definition~\ref{definition:Vogel-Cohn-Bousfield-local-derived-series}).
(In defining $\pi^{\{n\}}$ and $\pi^{(n)}$ we assume $\Gamma=\{e\}$ so
that the series are defined for any group~$\pi$.)  Let
$\{\pi^{(n)}_H\}$ be the Harvey derived
series~\cite{Cochran-Harvey:2004-1}.

\begin{theorem}
  \label{theorem:comparison}
  For any group $\pi$ and for any $n$, $ \pi^{n} \subset \pi^{\{n\}}
  \subset \pi^{(n)} \subset \pi_H^{(n)} $.
\end{theorem}

\begin{proof}
  We use an induction on~$n$.  For $n=0$, the conclusion holds since
  the initial (zeroth) terms of all the concerned series are $\pi$
  itself by definition.  Suppose the conclusion holds for~$n$.  Let
  $L(G)$ be the Cohn localization of the augmentation map $RG \to R$.
  For a group $G$ such that $\Z G$ is an Ore domain, we denote by
  $\K(G)$ the Ore localization $\Z G(\Z G-\{0\})^{-1}$ of~$\Z G$.  (In
  this proof $\K(G)$ is always well-defined whenever the notation
  $\K(G)$ is used.)

  We consider the following commutative diagram; note that, by
  definition, the $(n+1)$-st terms of the concerned series equal
  the kernels of the rows.  (To see that the kernel of the second row is
  $\pi^{\{n+1\}}=\Ker\{\pi \to \widehat\pi/\widehat\pi^{(n+1)}\}$, note
  that $\widehat\pi^{\{n\}} = \widehat\pi^{n}$.)
  \[
  \begin{diagram} \dgARROWLENGTH=1em \dgHORIZPAD=2mm \dgVERTPAD=1mm
    \def\hpi{\widehat\pi}
    \node{\pi^{n}} \arrow{e} \arrow{s,l}{\cap}
    \node{H_1\bigg(\pi;R\Big[\frac{\pi}{\pi^{n}}\Big]\bigg)
      \hbox to 0mm{$\displaystyle
        =\frac{\pi^{n}}{[\pi^{n},\pi^{n}]}\otimesover{\Z} R$
        \hss}
    }\arrow{s}
    \\
    \node{\pi^{\{n\}}} \arrow{e} \arrow[2]{s,l}{\cap}
    \node{H_1\bigg(\pi;R\Big[\frac{\pi}{\pi^{\{n\}}}\Big]\bigg)}
    \arrow[2]{s} \arrow{e}
    \node{H_1\bigg(\hpi;R\Big[\frac{\hpi}{\hpi^{\{n\}}}\Big]\bigg)
      \hbox to 0mm{$\displaystyle
        =\frac{\hpi^{\{n\}}}{[\hpi^{\{n\}},\hpi^{\{n\}}]}\otimesover{\Z} R$
        \hss}
    } \arrow{s}
    \\
    \node[3]{H_1\bigg(\hpi;R\Big[\frac{\hpi}{\hpi^{(n)}}\Big]\bigg)}
    \arrow{e}
    \node{H_1\bigg(\hpi;L\Big(\frac{\hpi}{\hpi^{(n)}}\Big)\bigg)}
    \\
    \node{\pi^{(n)}} \arrow{e} \arrow[2]{s,l}{\cap}
    \node{H_1\bigg(\pi;R\Big[\frac{\pi}{\pi^{(n)}}\Big]\bigg)}
    \arrow[2]{s} \arrow{e}
    \node{H_1\bigg(\pi;L\Big(\frac{\hpi}{\hpi^{(n)}}\Big)\bigg)}
    \arrow{ne,b}{\alpha} \arrow{s}
    \\
    \node[3]{H_1\bigg(\pi;L\Big(\frac{\hpi}{\hpi_H^{(n)}}\Big)\bigg)} \arrow{e}
    \node{H_1\bigg(\pi;\K\Big(\frac{\hpi}{\hpi_H^{(n)}}\Big)\bigg)}
    \\
    \node{\pi_H^{(n)}} \arrow{e}
    \node{H_1\bigg(\pi;\Q\Big[\frac{\pi}{\pi_H^{(n)}}\Big]\bigg)}\arrow{e}
    \node{H_1\bigg(\pi;\K\Big(\frac{\pi}{\pi_H^{(n)}}\Big)\bigg)}
    \arrow{ne,b}{\beta}
  \end{diagram}
  \]
  Note that in order to obtain the diagram we need the existence of
  the following natural maps:
  \begin{enumerate}
  \item $L(\widehat\pi/\widehat\pi_H^{(n)}) \to
    \K(\widehat\pi/\widehat\pi_H^{(n)})$.  Since
    $\widehat\pi/\widehat\pi_H^{(n)}$ is known to be a
    poly-torsion-free-abelian group~\cite{Cochran-Harvey:2004-1}, from
    Strebel's result \cite{Strebel:1974-1} it follows that the
    morphism $R[\widehat\pi/\widehat\pi_H^{(n)}] \to
    \K(\widehat\pi/\widehat\pi_H^{(n)})$ is an object of the category
    $\mathcal{C}^{\W_0}$ we described in
    Subsection~\ref{subsection:ring-localization} in order to define
    the Cohn localization.  Therefore, since
    $L(\widehat\pi/\widehat\pi_H^{(n)})$ is an initial object of
    $\mathcal{C}^{\W_0}$, the desired map exists.
  \item $\K(\pi/\pi_H^{(n)}) \to \K(\widehat\pi/\widehat\pi_H^{(n)})$.
    Since $\pi/\pi_H^{(n)} \to \widehat\pi/\widehat\pi_H^{(n)}$ is
    injective \cite{Cochran-Harvey:2004-1,Cha:2004-1}, it induces an
    injection $\K(\pi/\pi_H^{(n)}) \to
    \K(\widehat\pi/\widehat\pi_H^{(n)})$ between (skew) fields.
  \end{enumerate}
  Observe the following:
  \begin{enumerate}
  \item $\alpha$ is an isomorphism by
    Theorems~\ref{theorem:properties-of-group-localization}
    and~\ref{theorem:2-conn-prop-of-ring-localization}.
  \item $\beta$ is injective since any $\K(\pi/\pi_H^{(n)})$-module,
    especially $H_1(\pi;\K(\pi/\pi_H^{(n)}))$, is free and
    $\K(\widehat \pi/\widehat \pi_H^{(n)})$ is flat over its
    sub(skew)field $\K(\pi/\pi_H^{(n)})$.
  \end{enumerate}
  Using the above observations, the desired result for $n+1$ follows
  from a straightforward diagram chase.  The limit ordinal case is
  immediate.
\end{proof}

The next result deals with the special case of $R=\Z_p$.

\begin{theorem}
  \label{theorem:comparison-p-series}
  Suppose $R=\Z_p$ and $\pi$ is finitely generated.  Then $\pi^n =
  \pi^{\{n\}} = \pi^{(n)}$ and $\pi/\pi^n \cong
  \widehat\pi/\widehat\pi^n$ for $n<\infty$.
\end{theorem}

To prove this we need:

\begin{lemma}[Strebel, Levine]
  \label{lemma:Cohn-local-prop-of-p-group-ring}
  For any finite $p$-group $G$, the Cohn localization of the group
  ring $\Z_pG$ is $\Z_pG$ itself.
\end{lemma}

Lemma~\ref{lemma:Cohn-local-prop-of-p-group-ring} follows immediately
from Strebel's work \cite{Strebel:1974-1} or Levine's argument used in
the proof of \cite[Lemma 4.3, p89]{Levine:1994-1}.  We omit details.

\begin{proof}[Proof of Theorem~\ref{theorem:comparison-p-series}]
  We will show that $\pi^n = \pi^{(n)}$ and $\pi/\pi^n \cong
  \widehat\pi/\widehat\pi^n$ by an induction on~$n$.  For $n=0$,
  everything holds obviously.  Suppose the desired conclusion holds
  for~$n$.  We have
  \[
  \pi^{(n+1)} = \Ker\bigg\{\pi^n \to
  H_1\bigg(\pi;\Z_p\Big[\frac{\pi}{\pi^{n}}\Big]\bigg)
  \to   H_1\bigg(\pi;L\Big(\frac{\pi}{\pi^{n}}\Big)\bigg) \bigg \}.
  \]
  It is known that $\pi/\pi^n$ is a finite $p$-group if $\pi$ is
  finitely generated (e.g., see \cite{Cochran-Harvey:2007-01}.)  By
  Lemma~\ref{lemma:Cohn-local-prop-of-p-group-ring} stated above,
  $L(\pi/\pi^n)$ is equal to $\Z_p[\pi/\pi^n]$ and hence
  $H_1(\pi;L(\pi/\pi^{n})) = (\pi^n/[\pi^n,\pi^n]) \otimesover{\Z}
  \Z_p$.  It follows that $\pi^{(n+1)}$ is equal to~$\pi^n$.

  To prove that $\pi/\pi^{n+1} \cong \widehat\pi/\widehat\pi^{n+1}$,
  it suffices to show that $\pi^n/\pi^{n+1} \to
  \widehat\pi^n/\widehat\pi^{n+1}$ is an isomorphism, by applying the
  five lemma to the following diagram:
  \[
  \begin{diagram}\dgARROWLENGTH=1.5em
    \node{1} \arrow{e}
    \node{\pi^{n}/\pi^{n+1}} \arrow{e}\arrow{s}
    \node{\pi/\pi^{n+1}} \arrow{e}\arrow{s}
    \node{\pi/\pi^{n}} \arrow{e}\arrow{s}
    \node{1}
    \\
    \node{1} \arrow{e}
    \node{\widehat\pi^{n}/\widehat\pi^{n+1}} \arrow{e}
    \node{\widehat\pi/\widehat\pi^{n+1}} \arrow{e}
    \node{\widehat\pi/\widehat\pi^{n}} \arrow{e}
    \node{1}
  \end{diagram}
  \]
  Since $\pi^n \to H_1(\pi;L(\pi/\pi^n))$ is surjective,
  $\pi^n/\pi^{n+1} \cong H_1(\pi;L(\pi/\pi^n))$, and similarly
  $\widehat\pi^n/\widehat\pi^{n+1} \cong
  H_1(\widehat\pi;L(\pi/\pi^n))$.  The natural map
  $H_1(\pi;L(\pi/\pi^n)) \to H_1(\widehat\pi;L(\pi/\pi^n))$ is the map
  $\alpha$ in the proof of Theorem~\ref{theorem:comparison}, which was
  shown to be an isomorphism.  This completes the proof.
\end{proof}

\begin{remark}
  Using Theorem~\ref{theorem:comparison-p-series}, one can prove
  the result of Cochran and Harvey on the mod $p$ derived series
  \cite[Corollary 4.3]{Cochran-Harvey:2007-01} using the newer tools of this paper.
\end{remark}

\section{Examples of computation}
\label{section:examples}

We will discuss examples illustrating interesting properties and
elementary computational techniques of the local derived series.  Some
examples relate to later geometric applications.

The following lemma helps us compute lower terms of the
Vogel-Cohn $R$-local derived series.

\begin{lemma}
  \label{lemma:computation-of-lower-local-derived-series}
  Suppose $R$ is a subring of $\Q$ (and $\Gamma=\{e\}$).  Then lower
  successive quotients of the Vogel-Cohn local derived series
  $\{\pi^{(n)}\}$ of a group $\pi$ are given as follows:
  \[
  \frac{\pi^{(n)}}{\pi^{(n+1)}} =
  \begin{cases}
    H_1(\pi;\Z)/\text{$S$-torsion} & \text{for }n= 0 \\
    H_1\big(\pi;\Z[\pi/\pi^{(1)}]\big) / \text{$\Sigma$-torsion} &
    \text{for }n=1
  \end{cases}
  \]
  where $S=\{s\in \Z \mid s \text{ is invertible in } R\}$,
  $\Sigma=\{x\in \Z[\pi/\pi^{(1)}] \mid \epsilon(x)\in S\}$, and
  $\epsilon\colon \Z[\pi/\pi^{(1)}] \to \Z$ is the augmentation map.
\end{lemma}
Using the lemma, $\pi^{(n+1)}$ can be computed as the kernel of the
obvious map $\pi^{(n)} \to \pi^{(n)}/\pi^{(n+1)}$.  (For example, for a group $\pi$ with a presentation, the Reidemeister-Schreier
method may help in computing $\pi^{(n+1)}$ from
$\pi^{(n)}/\pi^{(n+1)}$.)

\begin{proof}
  For $n=0$, the conclusion is trivial.  For $n=1$, note that
  $\pi^{(1)}$ is the kernel of
  \[
  \pi \to H_1(\pi;\Z) \to H_1(\pi;R)=H_1(\pi;\Z)\otimes R.
  \]
  From this the conclusion follows, since the first map is surjective
  and the kernel of the second map is the $S$-torsion subgroup of
  $H_1(\pi;\Z)$.

  For $n=2$, first note that since $\widehat\pi/\widehat\pi^{(1)}$ is
  abelian, the Cohn localization $\Lambda$ of the augmentation map
  $R[\widehat\pi/\widehat\pi^{(1)}]\to R$ is isomorphic to
  $\widehat\Sigma^{-1} R[\widehat\pi/\widehat\pi^{(1)}]$ where
  $\widehat\Sigma = \{x\in \Z[\widehat\pi/\widehat\pi^{(1)}] \mid
  \widehat\epsilon(x) \in S\}$ and $\widehat\epsilon$ is the
  augmentation $\Z[\widehat\pi/\widehat\pi^{(1)}] \to \Z$.  Therefore
  $\pi^{(2)}$ is the kernel of
  \begin{multline}
    \pi^{(1)} \to H_1\Big(\pi;\Z\big[\frac{\pi}{\pi^{(1)}}\big]\Big)
    \xrightarrow{\alpha}
    H_1\Big(\pi;\Z\big[\frac{\widehat\pi}{\widehat\pi^{(1)}}\big]\Big)
    \\
    \xrightarrow{\beta} H_1(\pi;\Lambda) =
    H_1\Big(\pi;\Z\big[\frac{\widehat\pi}{\widehat\pi^{(1)}}\big]\Big)
    \otimesover{\Z[\frac{\widehat\pi}{\widehat\pi^{(1)}}]} \widehat
    \Sigma^{-1} \Z\big[\frac{\widehat\pi}{\widehat\pi^{(1)}}\big].
    \tag{$*$}
  \end{multline}
  For the last equality, note that in this case, $\Lambda$ a commutative localization, and therefor flat.

  Since $\pi/\pi^{(1)} \to \widehat\pi/\widehat\pi^{(1)}$ is
  injective, we view $\pi/\pi^{(1)}$ as a subgroup of
  $\widehat\pi/\widehat\pi^{(1)}$.  Let $W$ be a set of coset
  representatives.  As $\Z[\pi/\pi^{(1)}]$-modules, for each $g\in W$,
  $g\cdot \Z[\pi/\pi^{(1)}]$ is a submodule of
  $\Z[\widehat\pi/\widehat\pi^{(1)}]$ isomorphic to
  $\Z[\pi/\pi^{(1)}]$ and we have
  \begin{equation}
    \Z[\widehat\pi/\widehat\pi^{(1)}] = 
    \bigoplus_{g\in W} g\cdot \Z[\pi/\pi^{(1)}].
    \tag{$**$}
  \end{equation}
  Also, we have
  \begin{equation}
    H_1\Big(\pi;\Z\big[\frac{\widehat\pi}{\widehat\pi^{(1)}}\big]\Big)
    \cong H_1\Big(\pi;\Z\big[\frac{\pi}{\pi^{(1)}}\big]\Big)
    \otimesover{\Z[\frac{\pi}{\pi^{(1)}}]}
    \Z\big[\frac{\widehat\pi}{\widehat\pi^{(1)}}\big] \cong
    \bigoplus_{g\in W} g\cdot H_1\Big(\pi;\Z\big[\frac{\pi}{\pi^{(1)}}\big]\Big).
    \tag{$*{*}*$}
  \end{equation}
  In particular the map $\alpha$ in $(*)$ is injective and we can view
  $H_1(\pi;\Z[\pi/\pi^{(1)}])$ as a submodule of
  $H_1(\pi;\Z[\widehat\pi/\widehat\pi^{(1)}])$.

  We assert that for any $x\in H_1(\pi;\Z[\pi/\pi^{(1)}])$, $r\cdot x
  = 0$ for some $r\in \widehat\Sigma\subset
  \Z[\widehat\pi/\widehat\pi^{(1)}]$ if and only if $t\cdot x =0$ for
  some $t\in \Sigma \subset \Z[\pi/\pi^{(1)}]$.  The if part is clear.
  For the only if part, suppose $r\cdot x = 0$ for some $r\in
  \widehat\Sigma$.  We can write $r=\sum_{g\in W} g\cdot t_g$ with
  $t_g \in \Z[\pi/\pi^{(1)}]$ using~($**$).  Since $x\in
  H_1(\pi;\Z[\pi/\pi^{(1)}])$, summands of $r\cdot x=\sum_{g\in W}
  g\cdot (t_g \cdot x)$ are contained in distinct direct summands
  in~$(*{*}*)$.  It follows that $t_g \cdot x=0$ for each~$g$.  Let
  $t=\sum_{g\in W} t_g\in \Z[\pi/\pi^{(1)}]$.  Obviously $t\cdot x =
  0$, and $t$ is in $\widehat\Sigma$ since $\epsilon(t) = \sum
  \epsilon(t_i) = \widehat\epsilon(r)$.

  Note that for the map $\beta$ in $(*)$, $\Ker\beta$ is equal to the
  $\widehat\Sigma$-torsion submodule of
  $H_1(\pi;\Z[\widehat\pi/\widehat\pi^{(1)}])$.  Therefore, by the
  above assertion, $\Ker \beta\alpha =$ (pre-image of $\Ker \beta$) is
  exactly the $\Sigma$-torsion submodule of
  $H_1(\pi;\Z[\pi/\pi^{(1)}])$.  Since the first map in $(*)$ is
  surjective, from this it follows that $\pi^{(1)}/\pi^{(2)} \cong
  H_1(\pi;\Z[\pi/\pi^{(1)}])/\Sigma$-torsion.
\end{proof}

\begin{example}
  We will illustrate that in general the (integral) local derived
  series is distinct from the ordinary derived series and the Harvey
  series.  Let $\pi=\Z_p \rtimes \Z_2$, where $p>1$ is an odd integer
  and $\Z_p$ is viewed as a $\Z[\Z_2]$-module via the negation action
  of $\Z_2$ on~$\Z_p$.  That is, denoting by $t$ the generator of
  $\Z_2$, $t^{-1}r t = -r$ in $\pi$ for $r\in \Z_p$.

  First we will show that the ordinary derived series is as follows:
  \[
  \pi^n =
  \begin{cases}
    \pi & \text{ for }n=0, \\
    \Z_p & \text{ for }n=1, \\
    0 & \text{ for all }n\ge 2.
  \end{cases}
  \]
  
  \begin{proof}
    From the above definition of the $t$ action, it is easily seen
    that $\Z_p$ is in the kernel of the abelianization map and so
    $\pi/\pi^1 = \Z_2$, $\pi^1=\Z_p$.  Since $\pi^1$ is abelian,
    $\pi^n = 0$ for $n\ge 2$.
  \end{proof}
  
  On the other hand, for any finite group $\pi$ (including our case),
  it can be seen easily that the Harvey series $\{\pi_H^{(n)}\}$ is
  given as follows:
  \[
  \pi_H^{(n)} = \pi \quad\text{ for all } n.
  \]
  
  Now we will show that the Vogel-Cohn local derived series of $\pi$ (for
  $R=\Z$ and $\Gamma=\{e\}$) is given by
  \[
  \pi^{(n)} =
  \begin{cases}
    \pi & \text{ for }n=0, \\
    \Z_p & \text{ for }n\ge 1.
  \end{cases}
  \]
  \begin{proof}
    By Lemma~\ref{lemma:computation-of-lower-local-derived-series},
    $\pi^{(0)}=\pi$ and $\pi/\pi^{(1)} = H_1(\pi;\Z)=\Z_2$ and thus
    $\pi^{(1)} = \pi^1 = \Z_p$.  Also by
    Lemma~\ref{lemma:computation-of-lower-local-derived-series},
    $\pi^{(1)}/\pi^{(2)} = H_1(\pi;\Z[\Z_2])/\Sigma$-torsion where
    \[
    \Sigma=\{a+bt \in \Z[\Z_2] \mid \epsilon(a+bt)=a+b=1 \}.
    \]
    For $r\in \Z_p=\pi^1/\pi^2$ and $a+(1-a)t\in \Sigma$, the action
    is given by $(a+(1-a)t)\cdot r = (a-(1-a))r = (2a-1)r$.  In
    particular $r$ is annihilated by $a+(1-a)t\in \Sigma$ when
    $a=(p+1)/2$.  It follows that $H_1(\pi;\Z[\Z_2])$ is
    $\Sigma$-torsion.  Therefore $\pi^{(2)} = \pi^{(1)}$.
  \end{proof}
\end{example}  

\begin{example}
  We will illustrate that even the rational (i.e., the largest) local
  derived series may be strictly smaller than the Harvey series.  Our
  example is a group which is rationally Vogel-Cohn solvable but not
  Harvey solvable.  This results from using the Cohn localization to compute the  local derived series in place of the Ore localization used in the Harvey series, as our computation shows. Observe from this example that the Ore localization inverts excessively many ring elements, enlarging the size of the Harvey series.
  
  Let $\pi= A \rtimes \Z$, where
  \[
  A = \frac{\Z[t,t^{-1}]}{\langle(t-1)^2\rangle}
  \]
  and a preferred generator, say $t\in \Z$, acts on $A$ by
  multiplication.  Let $R$ be any subring of $\Q$ and $\pi^{(n)}$ be
  the $R$-coefficient Vogel-Cohn local derived series of~$\pi$.  We
  will show:
  \[
  \pi^{(n)} =
  \begin{cases}
    \pi & \text{ for }n=0, \\
    (t-1)A \cong \Z & \text{ for }n=1, \\
    0 & \text{ for }n\ge 2.
  \end{cases}
  \]

  \begin{proof}
    It can be seen that $[\pi,\pi] = (t-1)A \cong \Z$ and $H_1(\pi;\Z)
    =\pi/[\pi,\pi] \cong \Z\times \Z$ is a torsion free abelian group
    generated by (the cosets of) $1\in A$ and $t\in \Z$.  Therefore,
    by Lemma~\ref{lemma:computation-of-lower-local-derived-series}
    $\pi/\pi^{(1)} = H_1(\pi;\Z)=\Z\times \Z$ and so $\pi^{(1)} =
    [\pi,\pi]$.  Since $\pi^{(1)}$ is abelian,
    $\pi^{(1)}=H_1(\pi;\Z[\Z\times\Z])$.  It can be verified easily
    that the $\Z\times\Z$-action on $H_1(\pi;\Z[\Z\times\Z])\cong
    (t-1)A$ is trivial.  From this it follows that no nontrivial
    element in $H_1(\pi;\Z[\Z\times\Z])$ is $\Sigma$-torsion, where
    $\Sigma$ is as in
    Lemma~\ref{lemma:computation-of-lower-local-derived-series}.
    Therefore $0=\pi^{(2)}=\pi^{(3)}=\cdots$ by
    Lemma~\ref{lemma:computation-of-lower-local-derived-series}.
  \end{proof}
  
  Now we will compute the Harvey series to show
  \[
  \pi_H^{(n)} =
  \begin{cases}
    \pi & \text{ for }n=0, \\
    (t-1)A \cong \Z & \text{ for }n\ge 1.
  \end{cases}
  \]

  \begin{proof}
    Since $H_1(\pi;\Z)$ is torsion free, $\pi_H^{(1)} = [\pi,\pi]=(t-1)A
    \cong \Z$ and $\pi/\pi_H^{(1)}=H_1(\pi;\Z)=\Z\times\Z$.  To compute
    $\pi_H^{(2)}$, recall that $\pi_H^{(2)}$ is the kernel of
    \[
    \pi_H^{(1)} = H_1(\pi;\Z[\Z\times\Z]) \to H_1(\pi;\Z[\Z\times\Z])
    \otimesover{\Z[\Z\times\Z]} \K(\Z\times\Z)
    \]
    where $\K(\Z\times\Z)$ denotes the Ore localization of
    $\Z[\Z\times\Z]$.  Since the action of $\Z\times\Z$ is trivial on
    $\pi_H^{(1)}$, every element of $\pi_H^{(1)}$ is annihilated by any
    $r \ne 0 \in \Z[\Z\times\Z]$ which auguments to zero (i.e.,
    $\epsilon(r)=1$).  It follows that the above homomorphism is a zero
    map.  Therefore we have $\pi_H^{(1)} = \pi_H^{(2)} = \cdots$.
  \end{proof}
  
  We remark that the elements $r$ used above (such that
  $\epsilon(r)=0$) are inverted in the Ore localization but \emph{not}
  inverted in the Cohn localization.
\end{example}

\begin{example} For a finitely generated free group $F$, all the
  concerned commutator-type series agree; it is known that $F^n =
  F_H^{(n)}$ for any~$n$~\cite{Cochran-Harvey:2004-1}, and so by
  Theorem~\ref{theorem:comparison}, it follows that $F^n = F^{\{n\}} =
  F^{(n)} = F_H^{(n)}$ for any~$n$.
\end{example}

\begin{example}
  \label{example:spherical-space-form-groups}
  In what follows we describe the Vogel-Cohn local derived series of
  all the fundamental groups of 3-dimensional linear spherical space
  forms, that is, all finite groups acting on $S^3$ linearly without
  fixed points (see Milnor \cite{milnor:1957-2}).  Since
  Lemma~\ref{lemma:computation-of-lower-local-derived-series}
  (together with the Reidemeister-Schreier method if necessary)
  applies effectively, we omit details of the computation.  Let $R$ be
  a fixed subring of~$\Q$.

  \begin{enumerate}
 
  \item $\pi=Q_{8n} = \langle x,y \mid y^{2n} = x^2 = (xy)^2\rangle$
    with $n\ge 1$: write $n=2^r\cdot n_0$ with $n_0$ odd.  If $1/2
    \notin R$, then
    \[
    \frac{\pi^{(n)}}{\pi^{(n+1)}} \cong
    \begin{cases}
      \Z_2\times \Z_2 \text{ generated by $x$ and $y$}, &n=0, \\
      \Z_{2^{r+1}} \text{ generated by $y^2$}, &n=1,
    \end{cases}
    \]
    and
    \[
    \pi^{(n)} = 
    \begin{cases}
      \langle y^2 \mid (y^2)^{2n} = 1 \rangle \cong \Z_{2n}, &n=1, \\
      \langle y^{2^{r+2}} \mid (y^{2^{r+2}})^{n_0} = 1 \rangle \cong
      \Z_{n_0} \text{ generated by $y^{2^{r+2}}$}, &n=2.
    \end{cases}
    \]
    If $1/2 \in R$, then $\pi^{(n)}=\pi$ for all~$n$.

  \item $\pi = P_{48} = \langle x,y \mid x^2 = (xy)^3 = y^4$,
    $x^4=1\rangle$: if $1/2 \notin R$, then
    \[
    \frac{\pi^{(n)}}{\pi^{(n+1)}} \cong
    \begin{cases}
      \Z_2 \text{ generated by $x$}, &n=0,\\
      0 & n \ge 1,
    \end{cases}
    \]
    and 
    \[
    \pi^{(n)} =
      \langle u,v \mid u^2=v^3,\, uv^2uv^2u=v,\, u^4=1 \rangle
    \]
    for $n\ge 1$, where $u=y^2$, $v=xy$.  If $1/2 \in R$, then
    $\pi^{(n)}=\pi$ for all~$n$.

  \item $\pi = P_{120} = \langle x,y \mid x^2 = (xy)^3 = y^5$,
    $x^4=1\rangle$: $\pi$ is perfect and so $\pi^{(n)} = \pi$ for
    any~$n$ and any~$R$.

  \item $\pi = D_{2^k(2n+1)} = \langle x,y \mid x^{2k}=1,\,
    y^{2n+1}=1,\, xyx^{-1}=y^{-1} \rangle$: if $1/2 \notin R$, then 
    \[
    \frac{\pi^{(n)}}{\pi^{(n+1)}} \cong
    \begin{cases}
      \Z_{2^k} \text{ generated by $x$}, &n=0,\\
      0 & n \ge 1,
    \end{cases}
    \]
    and for $n\ge 1$,
    \[
    \pi^{(n)} = \langle y \mid y^{2n+1} = 1\rangle \cong \Z_{2n+1}.
    \]
    If $1/2\in R$, then $\pi^{(n)} = \pi$ for all~$n$.

  \item $\pi = P'_{8\cdot 3^k} = \langle x, y, z \mid
    x^2=(xy)^2=y^2,\, z^{-1}xz=y,\, z^{-1}yz=xy,\, z^{3^k}=1 \rangle$:
    if $1/3 \notin R$, then
    \[
    \frac{\pi^{(n)}}{\pi^{(n+1)}} \cong
    \begin{cases}
      \Z_{3^k} \text{ generated by $z$}, &n=0,\\
      0 & n \ge 1,
    \end{cases}
    \]
    and for $n\ge 1$,
    \[
    \pi^{(n)} = \langle x,y \mid xyx^{-1}=y^{-1},\, x^2 = y^2 \rangle.
    \]
    If $1/3\in R$, then $\pi^{(n)} = \pi$ for all~$n$.

  \item $\pi=\Z_d\times G$ where $G$ is either trivial or any group
    listed above: first consider $\pi=\Z_d$.  If $1/d\notin R$, then
    $\Z_d^{(n)}=0$ for any $n\ge 1$ .  If $1/d \in R$, then
    $\Z_d^{(n)}=\Z_d$ for any~$n$.  For the case of a product
    $\pi=\Z_d\times G$, the following lemma applies.
  \end{enumerate}
\end{example}

\begin{lemma}
  \label{lemma:series-for-direct-products}
  For the direct product $\pi=G\times H$ of two groups $G$ and $H$,
  $\pi^{(n)} = G^{(n)} \times H^{(n)}$ for any~$n$.
\end{lemma}

We remark that this product formula holds for the local derived series
determined by any homology localization triple.

\begin{proof}
  We use an induction on~$n$.  For $n=0$, the conclusion is trivial.
  The conclusion for a limit ordinal follows easily by taking limits.
  To consider the case of a successor ordinal, suppose $\pi^{(n)} =
  G^{(n)} \times H^{(n)}$.  Consider the following diagram:
  \[
  \begin{diagram} 
    \node{G^{(n)}\times H^{(n)}} \arrow{s,=}\arrow{e}
    \node{H_1\Big(G;L\Big(\dfrac{\widehat G}{\widehat G^{(n)}}\Big)\Big) \times
      H_1\Big(H;L\Big(\dfrac{\widehat H}{\widehat H^{(n)}}\Big)\Big)}
    \arrow{s,r}{\phi}
    \\
    \node{\pi^{(n)}} \arrow{e}
    \node{H_1\Big(\pi;L\Big(\frac{\widehat\pi}{\widehat\pi^{(n)}}\Big)\Big)}
  \end{diagram}
  \]
  where the map $\phi$ is the product of the maps induced by the
  inclusions of $G$ and $H$ into~$\pi$.  Note that there are induced
  maps on the homology coefficients by the functoriality of
  $(-)/(-)^{(n)}$ and~$L$.  By the functoriality of the group
  localization, the projection induces $\widehat\pi \to \widehat G$
  which is a left inverse of $\widehat G \to \widehat\pi$.  Therefore
  there is an induced map $\alpha\colon H_1(\pi;L(\widehat
  \pi/\widehat \pi^{(n)})) \to H_1(G;L(\widehat G/\widehat G^{(n)}))$.
  A similar argument applied to $H$ gives a map $\beta\colon
  H_1(\pi;L(\widehat \pi/\widehat \pi^{(n)})) \to H_1(H;L(\widehat
  H/\widehat H^{(n)}))$.  Then the map $(\alpha,\beta)$ on
  $H_1(\pi;L(\widehat \pi/\widehat \pi^{(n)}))$ is a left inverse of
  the map~$\phi$.  Therefore $\phi$ is injective.  From this it
  follows that $\pi^{(n+1)} = G^{(n+1)} \times H^{(n+1)}$ by a diagram
  chase.
\end{proof}

We finish this section by emphasizing that $\pi/\pi^{(n)}$ may have
nontrivial \emph{torsion} elements (e.g., see
Example~\ref{example:spherical-space-form-groups}).  In contrast to
this, the Harvey series quotient $\pi/\pi^{(n)}_H$ is always
torsion-free.  The existence of torsion elements enables us to make
some interesting applications of the local derived series to homology
cobordism of manifolds, as discussed in later sections.

\section{$\L2$-dimensional local property of von Neumann group
  algebras}
\label{section:$\L2$-dim-local-property}

In this section we will prove that the von Neumann group algebra of a group ring behaves much like a Cohn local module over the underlying group ring for an appropriate collection of groups.  This  ``$\L2$-dimensional Cohn local
property'' of the von Neumann group algebra $\N G$
(Theorem~\ref{theorem:$\L2$-dim-local}), has significant  geometric applications
related to homological properties.  The essential new contribution is
that it holds even for a class of groups $G$ that may have
\emph{torsion} (as well as infinite order elements).  We start by
defining the notion of an $\L2$-equivalence.

\subsection{$\L2$-dimension and $\L2$-equivalences}

Suppose $G$ is a countable group with von Neumann group algebra~$\N
G$.  Following L\"uck's book \cite[Chapter~6]{Lueck:2002-1}, the
$\L2$-dimension function
\[
\ldim\colon \{ \text{$\N G$-modules} \} \to [0,\infty] =
[0,\infty) \cup \{\infty\}
\]
is defined as follows.  For a finitely generated projective $\N
G$-module $P$, $\ldim P$ is defined to be the von Neumann trace of a
square matrix over $\N G$ whose row space is isomorphic to~$P$.  For
an arbitrary module $M$ over $\N G$, $\ldim M=\sup\{\ldim P \mid P$ is
a finitely generated projective submodule of $M\}$.

In this paper we will frequently use the following basic properties of
the $\L2$-dimension function:
\begin{enumerate}
\item $\ldim \N G=1$ and $\ldim 0 = 0$.
\item If $0\to M'\to M \to M'' \to 0$ is exact, then $\ldim M = \ldim
  M' + \ldim M''$.  In particular, if $N$ is either a submodule or a
  homomorphic image of $M$, then $\ldim N \le \ldim M$.
\end{enumerate}
For proofs, see \cite[Chapter~6]{Lueck:2002-1}.  Here we adopt the
usual convention $\infty+d=\infty$ for any $d\in [0,\infty]$.

\begin{definition}
  An $\N G$-module homomorphism $f\colon M\to N$ is said to be an
  \emph{$\L2$-equivalence} if $\ldim \Ker f = 0 = \ldim \Coker f$.
\end{definition}

\begin{remark}
  \leavevmode\Nopagebreak
  \begin{enumerate}
  \item If $f\colon M \to N$ is an $\L2$-equivalence, $\ldim M = \ldim
    N= \ldim \Im f$.
  \item If $\ldim M = \ldim N < \infty$, then we have $\ldim \Ker f =
    \ldim \Coker f$ for any $f\colon M \to N$.  In particular in this
    case $f$ is an $\L2$-equivalence if and only if either $\Ker f$ or
    $\Coker f$ has $\ldim=0$.
  \end{enumerate}
\end{remark}

A group $G$ is called \emph{amenable} if $\N G \otimesover{\C G} \C$
is nonzero.  For other equivalent definitions (which are more often
used in the literature) and related discussions, the reader is
referred to Paterson's book \cite{Paterson:1988-1} and L\"uck's book
\cite[p.~256]{Lueck:2002-1}.  For our purpose, the following fact on
amenable groups will play an essential role.

\begin{theorem}[{L\"uck~\cite[Theorem~6.37]{Lueck:2002-1}}]
  \label{theorem:NG-is-$\L2$-dim-flat}
  Suppose $G$ is an amenable group.  Then $\N G$ is
  $\L2$-dimension flat over $\C G$, that is, $\ldim
  \Tor_p^{\C G}(M,\N G)=0$ for any $\C G$-module~$M$ and
  for any $p>0$.
\end{theorem}

\begin{lemma}
  \label{lemma:injection-induces-$\L2$-equiv}
  Suppose $G$ is amenable.  If $M$ and $N$ are free $\C G$-modules
  with the same finite rank and $f\colon M \to N$ is an injection,
  then $f\otimes 1_{\N G} \colon M\otimes_{\C G} \N G \to N\otimes_{\C
    G} \N G$ is an $\L2$-equivalence.
\end{lemma}

\begin{proof}
  Let $C=\Coker f$.  From $0 \to M\to N \to C \to 0$ we obtain
  \[
  0 \to \Tor^{\C G}_1 (C,\N G) \to M\otimesover{\C G}
  \N G \to N\otimesover{\C G} \N G \to C
  \otimesover{\C G} \N G \to 0.
  \]
  By Theorem~\ref{theorem:NG-is-$\L2$-dim-flat}, $\ldim \Tor^{\C G}_1
  (C,\N G)=0$.  Since $M\otimes_{\C G} \N G$ and
  $N\otimes_{\C G} \N G$ are finitely generated free
  $\N G$-modules with the same rank, they have the same finite
  $\ldim$.  It follows that $f\otimes 1_{\N G}$ is an
  $\L2$-equivalence.
\end{proof}

The following is an $\L2$-equivalence analogue of the observation that
a chain map which is an isomorphism induces an isomorphism on the
homology.

\begin{lemma}
  \label{lemma:homology-$\L2$-equiv}
  Suppose $\phi_*\colon C_* \to C_*$ is a chain map of a chain complex
  $C_*$ over~$\N G$.  If $\ldim C_i < \infty$ and $\phi_i\colon
  C_i \to C_i$ is an $\L2$-equivalence for some~$i$, then the induced
  map $(\phi_i)_* \colon H_i(C_*) \to H_i(C_*)$ is an
  $\L2$-equivalence.
\end{lemma}

\begin{proof}
  Let $Z_i, B_i \subset C_i$ be the submodules of cycles and
  boundaries, respectively.  Consider $\phi_i|_{Z_i}\colon Z_i\to
  Z_i$.  Since $\Ker\phi_i|_{Z_i} = Z_i \cap \Ker \phi_i$ is a
  submodule of $\Ker \phi_i$ and $\ldim \Ker \phi_i = 0$, we have
  $\ldim \Ker\phi_i|_{Z_i} = 0$.  Since $\ldim Z_i \le \ldim C_i
  <\infty$, it follows that $\phi_i|_{Z_i}$ is an $\L2$-equivalence.
  \[
  \Coker\{ (\phi_i)_*\colon Z_i/B_i \to Z_i/B_i \} =
  Z_i/(B_i+\phi_i(Z_i))
  \]
  is a homeomorphic image of $\Coker \phi_i|_{Z_i} = Z_i/\phi_i(Z_i)$
  and therefore
  \[
  \ldim \Coker (\phi_i)_* \le \ldim \Coker
  \phi_i|_{Z_i} = 0.
  \]
  Since $\ldim Z_i/B_i \le \ldim C_i <\infty$, it follows that
  $(\phi_i)_*$ is an $\L2$-equivalence.
\end{proof}

\subsection{$\N G$-homology and Strebel's class $D(R)$}

Suppose $R$ is a commutative ring.  We always assume that a
commutative ring $R$ has unity and $1\ne 0$ in $R$, i.e., the natural
map $\Z \to R$ is nonzero.  In \cite{Strebel:1974-1}, Strebel defined
and studied the class $D(R)$ of groups $G$ with the following
property: whenever $f\colon M\to N$ is a homomorphism between
projective $RG$-modules such that $f\otimes 1_R\colon M\otimesover{RG}
R\to N \otimesover{RG} R$ is injective, $f$ itself is injective.
(Here $R$ is viewed as an $RG$-module with trivial $G$-action.)


In the following theorem, we relate the class $D(R)$ with the
$\L2$-dimension of $\N G$-homology modules.

\begin{theorem}
  \label{theorem:$\L2$-dim-local}
  Suppose $R$ is a commutative ring, $G$ is an amenable group, and
  $G \to \Gamma$ is a group homomorphism with kernel in~$D(R)$.
  Suppose $C_*$ is a bounded below chain complex over $\Z G$ such that
  $C_i$ is finitely generated and free as a $\Z G$-module for $i\le n$.  If
  $H_i(C_*\otimesover{\Z G} R\Gamma)=0$ for $i\le n$, then $\ldim
  H_i(C_*\otimesover{\Z G} \N G) = 0$ for $i\le n$.
\end{theorem}

\begin{remark}
  \leavevmode\Nopagebreak
  \begin{enumerate}
  \item Theorem~\ref{theorem:$\L2$-dim-local} will be useful in proving
    our results on $R\Gamma$-coefficient homology cobordism.
    Our primary examples applying this theorem will use the rings $R =
    \Z_p$ and $R =\Q$ in Theorem~\ref{theorem:$\L2$-dim-local}.  One
    can consider subrings $R \subset\Q$ (including $\Z$ and
    $\Z_{(p)}$), but this serves no purpose since in this case,
    $D(R)=D(\Q)$.  On the other hand, the classes $D(\Q)$ and
    $D(\Z_p)$ are distinct.  For example, a finite $p$-group is in
    $D(\Z_p)$ but not in $D(\Q)$.  See also
    Lemma~\ref{lemma:amenable-and-D(R)}.
  \item It is well known that a similar result to
    Theorem~\ref{theorem:$\L2$-dim-local} holds for the
    \emph{Cohn localization} of~$\Z G
    \to \Z\Gamma$. (See, for instance, \cite{Vogel:1982-1}.) In this sense, we can interpret the conclusion of
    Theorem~\ref{theorem:$\L2$-dim-local} as follows: the von Neumann
    group ring $\N G$ is an \emph{$\L2$-dimension Cohn local $\Z
      G$-module}.
  \end{enumerate}
\end{remark}

\begin{proof}[Proof of Theorem~\ref{theorem:$\L2$-dim-local}]
  We will show that the zero map of $H_i(C_*\otimesover{\Z G}
  \N G)$ into itself is an $\L2$-equivalence.  From this the
  desired conclusion follows immediately.

  Our argument is similar to Vogel's argument used in
  \cite{Vogel:1982-1} but requires
  Lemma~\ref{lemma:injection-induces-$\L2$-equiv} at a crucial point.
  Since $H_i(C_*\otimesover{\Z G} R\Gamma)=0$ for $i\le n$, there is a
  partial chain homotopy
  \[
  s_i\colon C_i\otimesover{\Z G}R\Gamma
  \to C_{i+1}\otimesover{\Z G}R\Gamma \quad (i\le n)
  \]
  such that $\partial s+s\partial = 1_{C_*\otimesover{\Z G} R\Gamma}$.
  It is easily seen that $s_i$ lifts to a partial chain homotopy
  $D_i\colon C_i \to C_{i+1}$ $(i\le n)$ such that $D\otimesover{\Z G}
  1_{R\Gamma} = s$.  Let $N=\Ker\{G \to \Gamma\}$.  For the partial
  chain map $\phi = \partial D + D\partial$ on $C_*$,
  \[
  (\phi_i \otimesover{\Z} 1_R) \otimesover{R N} 1_R\colon
  (C_i\otimesover{\Z} R)\otimesover{R N} R \to (C_i\otimesover{\Z}
  R)\otimesover{R N} R
  \]
  can be identified with
  \[
  \phi_i \otimesover{\Z G} 1_{R\Gamma}\colon C_i\otimesover{\Z G} R\Gamma \to
  C_i\otimesover{\Z G} R\Gamma
  \]
  which is the identity map ($i\le n$).  Therefore, since $N$ is in
  $D(R)$,
  \[
  \phi_i \otimesover{\Z} 1_{R} \colon C_i \otimesover{\Z} R \to C_i
  \otimesover{\Z} R
  \]
  is injective.  Note that for any (possibly infinite) index set $A$,
  a homomorphism $\Z^A \to \Z^A$ is injective if the induced map $R^A
  \to R^A$ is injective.  From this it follows that $\phi_i \colon C_i
  \to C_i$ is injective.  Now, by
  Lemma~\ref{lemma:injection-induces-$\L2$-equiv},
  \[
  \phi_i\otimes 1_{\N G} \colon C_i\otimesover{\Z G} \N G
  \to C_i\otimesover{\Z G} \N G
  \]
  is an $\L2$-equivalence.  Therefore, by
  Lemma~\ref{lemma:homology-$\L2$-equiv} it follows that $\phi_i$ induces
  an $\L2$-equivalence
  \[
  \phi_*\colon H_i(C_*\otimesover{\Z G} \N G) \to
  H_i(C_*\otimesover{\Z G} \N G).
  \]
  Because of $D_*$, the induced map $\phi_*$ is zero for $i\le n$.
\end{proof}

Note that if $\Gamma$ is trivial in
Theorem~\ref{theorem:$\L2$-dim-local}, or more generally if $\Gamma$
is amenable, then $G$ is amenable if and only if $\Ker\{G\to \Gamma\}$
is amenable, since the class of amenable groups is closed under
extensions and taking subgroups and quotients.  Therefore the
hypothesis of Theorem~\ref{theorem:$\L2$-dim-local} leads us to
consider the class of groups which are amenable and in~$D(R)$.  We
list some known cases (e.g., see \cite{Paterson:1988-1},
\cite{Lueck:2002-1}, and \cite{Strebel:1974-1}).

\begin{enumerate}
\item Poly-torsion-free-abelian (PTFA) groups are amenable and in
  $D(R)$ for any ring~$R$.
\item Finite $p$-groups are amenable and in~$D(\Z_p)$.
\item A direct limit of amenable groups (resp.\ groups in $D(R)$) is
  amenable (resp.\ in~$D(R)$).
\item If $G$ admits a subnormal series
  \[
  G = G_0 \supset G_1 \supset \cdots \supset G_n =
  \{e\}
  \]
  such that each $G_i/G_{i+1}$ is amenable (resp.\ in $D(R)$), then
  $G$ is amenable (resp.\ in~$D(R)$).
\end{enumerate}

The following consequence will be used later.  We say that a group
\emph{$G$ has no torsion coprime to $p$} if the order of any finite
order element in $G$ is a power of~$p$.  (Note that $G$ may have
infinite order elements and may be infinitely generated.)

\begin{lemma}
  \label{lemma:amenable-and-D(R)}
  Suppose $G$ is a group admitting a subnormal series
  \[
  G = G_0 \supset G_1 \supset \cdots \supset G_n =
  \{e\}
  \]
  whose quotients $G_i/G_{i+1}$ are abelian.
  \begin{enumerate}
  \item If every $G_i/G_{i+1}$ has no torsion coprime to $p$, then $G$ is
    amenable and in $D(\Z_p)$.
  \item If every $G_i/G_{i+1}$ is torsion free, then $G$ is amenable
    and in $D(R)$ for any commutative ring~$R$.
  \end{enumerate}
\end{lemma}

\begin{proof}
  (1) It suffices to show that $G_i/G_{i+1}$ is amenable and
  in~$D(R)$.  Note that $G_i/G_{i+1}$ is the direct limit of its
  finitely generated subgroups.  Since $G_i/G_{i+1}$ is abelian and
  has no torsion coprime to $p$, each finitely generated subgroup $H$
  of $G_i/G_{i+1}$ is a direct sum of abelian $p$-groups and free
  abelian groups.  Therefore $H$ is amenable and in $D(\Z_p)$.  It
  follows that $G_i/G_{i+1}$ is amenable and in $D(\Z_p)$.
  (Amenability of $G_i/G_{i+1}$ also follows immediately from that it
  is abelian.)
   


  (2) In this case, $G_i/G_{i+1}$ is torsion-free and so is a direct
  limit of finitely generated free abelian groups.  It follows that
  $G$ is amenable and in $D(R)$.
\end{proof}

\section{$\L2$-invariants and homology cobordism}
\label{section:$\L2$-inv-homology-cobordism}

For a CW-complex $X$ and a group homomorphism $\phi\colon \pi_1(X)\to
G$, we define the \emph{$\L2$-Betti number} by $b^{(2)}_i(X,\phi) =
\ldim H_i(X;\N G)$.  When $X=M$ is a $(4k-1)$-dimensional manifold, we
denote by $\rho^{(2)}(M,\phi)$ the \emph{$\L2$-signature defect} of a
bounding $4k$-manifold over~$G$.

For the reader's convenience, we sketch the definition of
$\rho^{(2)}(M,\phi)$.  It is known that there exist a compact
$4k$-manifold $W$ with $\partial W=M$, an \emph{injection} $G \to H$
of groups, and a homomorphism $\psi\colon \pi_1(W) \to H$ such that
the following diagram commute:
\[
\begin{diagram}
  \node{\pi_1(M)} \arrow{e,t}{\phi} \arrow{s,l}{i_*}
  \node{G} \arrow{s,..}
  \\
  \node{\pi_1(W)}\arrow{e,b,..}{\psi}
  \node{H}
\end{diagram}
\]
Applying the spectral theory to the $\N H$-coefficient intersection
form
\[
H_{2k}(W;\N H)\times H_{2k}(W;\N H) \to \N H
\]
$H_{2k}(W;\N H)$ can be written as an orthogonal sum $V_+ \oplus V_-
\oplus V_0$ where the intersection form on $V_+$, $V_-$, and $V_0$ are
positive definite, negative definite, and zero, respectively.  The
\emph{$\L2$-signature} of $(W,\psi)$ is defined to be
\[
\lsign(W,\psi)=\ldim V_+ - \ldim V_-.
\]
Now the $\rho$-invariant $\rho^{(2)}(M,\phi)$ is defined to be the
\emph{$\L2$-signature defect} of $W$, namely
\[
\rho^{(2)}(M,\phi)=\lsign(W,\psi)-\sign(W)
\]
where $\sign(W)$ is the ordinary signature of~$W$.  It is well known
that $\rho^{(2)}(M,\phi)$ is a well-defined real-valued invariant of
$(M,\phi)$, which agrees with the Cheeger-Gromov invariant.  For more
details and related discussions, the reader may be referred to, for
example, Cochran-Orr-Teichner~\cite{Cochran-Orr-Teichner:1999-1},
Chang-Weinberger~\cite{Chang-Weinberger:2003-1},
L\"uck~\cite{Lueck:2002-1}, Harvey~\cite{Harvey:2006-1}, and
Cha~\cite{Cha:2006-1}.

We also need the induction property: if $\phi\colon \pi_1(M) \to
G$ is a homomorphism and $f \colon G \to H$ is an injection,
then $\rho^{(2)}(M,\phi) = \rho^{(2)}(M, f\phi)$.

We say that a (connected) complex $X$ is over a group $\Gamma$ if $X$
is endowed with a homomorphism $\pi_1(X) \to \Gamma$.  For two closed
manifolds $M$ and $M'$ over $\Gamma$ and a commutative ring $R$, a
bordism $W$ over $\Gamma$ between $M$ and $M'$ is called an
\emph{$R\Gamma$-homology cobordism} if the inclusions of $M$ and $M'$
into $W$ induce isomorphisms on $H_*(-;R\Gamma)$.  If such $W$ exists,
$M$ and $M'$ are said to be \emph{$R\Gamma$-homology cobordant}.

\begin{theorem}
  \label{theorem:homology-invariance-$\L2$-sign}
  Suppose $R$ is a commutative ring, $G$ is an amenable group, and $G
  \to \Gamma$ is a homomorphism with kernel in~$D(R)$.  Suppose $W$ is
  an $R\Gamma$-homology cobordism between closed manifolds $M$ and
  $M'$ over~$\Gamma$.  If $\phi\colon \pi_1(M) \to G$ and $\phi'\colon
  \pi_1(M') \to G$ are restrictions of a homomorphism $\psi\colon
  \pi_1(W) \to G$, then the following hold:
  \begin{enumerate}
  \item $b^{(2)}_i(M,\phi) = b^{(2)}_i(M',\phi')$ for all~$i$.
  \item When $M$ and $M'$ are $(4k-1)$-dimensional,
    $\rho^{(2)}(M,\phi) = \rho^{(2)}(M',\phi')$.
  \end{enumerate}
\end{theorem}

Note that $W$ is a $\Z_p\Gamma$-homology cobordism if and only if $W$
is a $\Z_{(p)}\Gamma$-homology cobordism.  Using this fact we can
combine the $\Z_{(p)}$-coefficient Vogel-Cohn local derived series
with the $\L2$-signatures associated to amenable groups in $D(\Z_p)$.
A precise statement is as follows.

\begin{theorem}
  \label{theorem:homology-invariance-der-series}
  Let $R$ be either $\Z_p$, $\Z_{(p)}$, or $\Q$.  For a closed
  manifold $M$ over an amenable group $\Gamma$, view
  $\pi = \pi_1(M)$ as a group over $\Gamma$ and denote the
  $R$-coefficient Vogel-Cohn local derived series of $\pi$ over
  $\Gamma$ by $\{\pi^{(n)}\}$.  For the canonical map $\phi_n\colon
  \pi \to \pi / \pi^{(n)}$, $b^{(2)}(M,\phi_n)$ and
  $\rho^{(2)}(M,\phi_n)$ are $R\Gamma$-homology cobordism invariants
  of~$M$ for any~$n<\infty$.  In particular, when $\Gamma$ is trivial,
  $b^{(2)}(M,\phi_n)$ and $\rho^{(2)}(M,\phi_n)$ are always
  $R$-homology cobordism invariants.
\end{theorem}

We remark that in many cases the group $\pi_1(M) / \pi_1(M)^{(n)}$ may
have torsion elements.  (For instance, see our computation for
spherical space forms:
Example~\ref{example:spherical-space-form-groups}.)  Prior to this work, only
$\L2$-signatures associated to \emph{poly-torsion-free-abelian} groups
have been known to be useful in studying homology cobordism.
Theorem~\ref{theorem:homology-invariance-der-series} is the first
result on the homology cobordism invariance of $\L2$-signatures
associated to \emph{non-torsion-free} groups.

\begin{proof}[Proof of Theorem~\ref{theorem:homology-invariance-$\L2$-sign}]
  Since $H_i(W,M;R\Gamma)=0$, applying
  Theorem~\ref{theorem:$\L2$-dim-local} to the chain complex
  $C_*=C_*(W,M;\Z G)$ we obtain $\ldim H_i(W,M;\N G)=0$ for all~$i$.
  From the $\N G$-coefficient homology long exact sequence for
  $(W,M)$, it follows that the map $H_i(M;\N G) \to H_i(W;\N G)$ is an
  $\L2$-equivalence.  The same argument applies to $M'$.  Therefore
  $b^{(2)}_i(M,\phi) = b^{(2)}_i(W,\psi) = b^{(2)}_i(M',\phi')$.

  For (2), note that $\rho^{(2)}(M,\phi) - \rho^{(2)}(M',\phi')$ is
  the difference between the $\L2$-signature of the $\N G$-intersection
  form of $W$ and the signature of the ordinary (untwisted)
  intersection form of~$W$.  We will show that both signatures are
  zero.  
  
  The non-singular part of the $\N G$-intersection form is
  supported by the cokernel of $H_{2k}(M;\N G) \to H_{2k}(W;\N G)$
  which is contained in $H_{2k}(W,M;\N G)$.  Since $\ldim
  H_{2k}(W,M;\N G)$ is zero, the $\L2$-signature of $W$ is zero.
  Similarly, the ordinary signature of $W$ is zero, since
  $H_{2k}(W,M;\Q)=0$.
\end{proof}

\begin{proof}[Proof of
  Theorem~\ref{theorem:homology-invariance-der-series}] We first
  consider the case $R=\Z_{(p)}$.  Suppose $W$ is a $R\Gamma$-homology
  cobordism between $M$ and~$M'$.  Let $G=\pi_1(W)/\pi_1(W)^{(n)}$,
  and let $\phi\colon \pi_1(M) \to \pi_1(W) \to G$ and $\phi'\colon
  \pi_1(M') \to \pi_1(W) \to G$ be the compositions.  By the
  injectivity theorem, $\pi_1(M)/\pi_1(M)^{(n)}$ injects into $G$.
  Therefore by the induction property, $\rho^{(2)}(M,\phi_n)$ is equal
  to $\rho^{(2)}(M,\phi)$, and similarly for~$M'$.  So it suffices to
  show that $\rho^{(2)}(M,\phi) = \rho^{(2)}(M',\phi')$.

  By the following lemma, $\Ker\{G\to \Gamma\}=
  \pi_1(W)^{(0)}/\pi_1(W)^{(n)}$ is amenable and in~$D(\Z_p)$.  From
  the amenability of $\Ker\{G\to \Gamma\}$ and~$\Gamma$, the
  amenability of $G$ follows.  Therefore by applying
  Theorem~\ref{theorem:homology-invariance-$\L2$-sign}, we obtain
  $\rho^{(2)}(M,\phi) = \rho^{(2)}(M',\phi')$.

  It can be seen that a similar argument works for $R=\Q$ and~$\Z_p$.
\end{proof}

\begin{lemma}
  Suppose $\pi$ is a group over~$\Gamma$.
  \begin{enumerate}
  \item Let $R$ be either $\Z_p$ or $\Z_{(p)}$ and $\pi^{(n)}$ be the
    $R$-coefficient Vogel-Cohn local derived series of $\pi$ over
    $\Gamma$.  Then $\pi^{(0)}/\pi^{(n)}$ is amenable and in $D(\Z_p)$
    for any~$n <\infty$.
  \item Let $\pi^{(n)}$ be the rational Vogel-Cohn local derived
    series of $\pi$ over $\Gamma$ and $R$ be any commutative ring.
    Then $\pi^{(0)}/\pi^{(n)}$ is amenable and in $D(R)$ for any~$n
    <\infty$.
  \end{enumerate}
\end{lemma}

\begin{proof}
  For (1), consider the normal series $\{\pi^{(i)}/\pi^{(n)}\}_{0\le
    i\le n}$ of~$\pi^{(0)}/\pi^{(n)}$.  The $i$-th quotient is
  isomorphic to $\pi^{(i)}/\pi^{(i+1)}$, which is the image of a map
  factoring through $H_1(\pi;R[\pi/\pi^{(i)}])$ by the definition
  of~$\pi^{(i+1)}$.  Since $H_1(\pi;R[\pi/\pi^{(i)}])$ has no torsion
  coprime to $p$ for both $R=\Z_p$ and $\Z_{(p)}$, so does
  $\pi^{(i)}/\pi^{(i+1)}$.  It follows that $\pi^{(0)}/\pi^{(n)}$ is
  amenable and in $D(\Z_p)$ by Lemma~\ref{lemma:amenable-and-D(R)}.
  The second part is proved similarly.
\end{proof}

\section{Applications}
\label{section:applications}

\subsection{Distinct homology cobordism types with the same simple
  homotopy type}
\label{subsection:distinct-homology-cobordism-types}

In~\cite{Chang-Weinberger:2003-1}, Chang and Weinberger proved the
following result on homeomorphism types of manifolds having the same
simple homotopy type using $\L2$-signature invariants:

\begin{theorem}[Chang-Weinberger~\cite{Chang-Weinberger:2003-1}]
  Suppose $M$ is a closed $(4k-1)$-manifold with $\pi=\pi_1(M)$, $k\ge
  2$.  If $\pi$ has a nontrivial torsion element, then there
  exist infinitely many closed $(4k-1)$-manifolds $M_0=M$, $M_1$,
  $M_2,\ldots$ such that each $M_i$ is simple homotopy equivalent and
  tangentially equivalent to $M$ but $M_i$ and $M_j$ are not
  homeomorphic for any $i\ne j$.
\end{theorem}

Using our results, we will prove a homology cobordism version of the
Chang-Weinberger theorem.  

\begin{theorem}
  \label{theorem:homology-cobordism-of-homotopy-equivalent-manifolds}
  Suppose $M$ is a closed $(4k-1)$-manifold with $\pi=\pi_1(M)$, $k\ge
  2$.  Let $p$ be prime and $\pi^{(n)}$ be the $\Z_p$ or
  $\Z_{(p)}$-coefficient Vogel-Cohn local derived series of $\pi$.  If
  $\pi$ has a torsion element which remains nontrivial in
    $\pi/\pi^{(n)}$ for some $n$, then there exist infinitely many
  closed $(4k-1)$-manifolds $M_0=M$, $M_1$, $M_2,\ldots$ such that
  each $M_i$ is simple homotopy equivalent and tangentially equivalent
  to $M$ but $M_i$ and $M_j$ are not homology cobordant for any $i\ne
  j$.
\end{theorem}

Note that $\pi/\pi^{(n)}$ has no torsion coprime to $p$ but may have
$p$-torsion elements.  We also remark that the $R\Gamma$-homology
generalizations of
Theorem~\ref{theorem:homology-cobordism-of-homotopy-equivalent-manifolds}
are true for $R=\Z$, $\Z_p$, and $\Z_{(p)}$ when $M$ is over an
amenable group~$\Gamma$.  (Indeed the same proof works.)

In the $\Z_p$-coefficient case of
Theorem~\ref{theorem:homology-cobordism-of-homotopy-equivalent-manifolds},
$\pi^n=\pi^{(n)}$ by Theorem~\ref{theorem:comparison-p-series}.
Therefore
Theorem~\ref{theorem:homology-cobordism-of-homotopy-equivalent-manifolds}
specializes to the following result, which can also be proved directly
by using the homology cobordism invariance of $\L2$-signatures
associated to finite $p$-groups:

\begin{corollary}
  \label{corollary:homology-cobordism-of-homotopy-equivalent-manifolds-mod-p}
  Suppose $M$ is a closed $(4k-1)$-manifold, $k>2$, and let
  $\{\pi^n\}$ be the $\Z_p$-derived series of $\pi=\pi_1(M)$.  If
    $\pi$ has a torsion element which is not contained in $\pi^n$ for
    some $n$, then the conclusion of
  Theorem~\ref{theorem:homology-cobordism-of-homotopy-equivalent-manifolds}
  holds.
\end{corollary}

\begin{proof}
  [Proof of
  Theorem~\ref{theorem:homology-cobordism-of-homotopy-equivalent-manifolds}]
  Following the argument of the proof of
  \cite[Theorem~1]{Chang-Weinberger:2003-1}, we consider the action of
  the surgery obstruction group $L_{4k}(\pi)$ on the structure set
  $\mathcal{S}(M)$ of simple homotopy equivalences of closed
  $(4k-1)$-manifolds into~$M$: for each $a\in L_{4k}(\pi)$, by the
  Wall realization theorem there is a $4k$-dimensional bordism, say
  $W_a$, between $M$ and another element $\mathcal{S}(M)$, say
  (represented by) $M_a$, over $\pi$ such that $a$ is represented by
  the $\Z\pi$-coefficient intersection form of~$W_a$.  We will show
  that we can choose infinitely many $a$, including $a=0$ (i.e.,
  $M=M_a$ for some $a$) so that $M_a$ is not homology cobordant to
  $M_{a'}$ for any $a\ne a'$.

  Let $G=\pi/\pi^{(n)}$ and $\phi\colon \pi_1(M) \to G$ and
  $\phi_a\colon \pi_1(M_a) \to G$ be the obvious maps.  Then
  $\rho^{(2)}(M,\phi)-\rho^{(2)}(M_a,\phi_a)$ is equal to the
  $\L2$-signature defect of the intersection form on $H_{2k}(W_a;\N
  G)$.  The proof of \cite[Theorem~1]{Chang-Weinberger:2003-1} shows
  the following: for any group $G$ which is not torsion free, there
  exist infinitely many $o\in L_{4k}(G)$ such that the $\L2$-signature
  defects of the forms representing $o$ are all nonzero and distinct.
    The elements $o$ are constructed as follows.  Switching to the
    $\Q$-coefficients and considering forms over projective modules
    makes no difference, since the map of $L_{4k}(G)=L_{4k}(\Z G)$
    into the relevant $L$-group $L_{4k}(\Q G)$ induces an isomorphism
    modulo 8-torsion for any~$G$.  (We omit decorations in the
    $L$-group notation.)  An order $r>1$ element in $G$ determines a
    map $L_{4k}(\Q[\Z_r]) \to L_{4k}(\Q G)$, and the images of direct
    sums of the $1\times 1$ form $[1]$ on the projective
    $\Q[\Z_r]$-module~$\Q$ has the desired property.

    Returning to our case, we have an element in $\pi$ with finite
    order, say $d$, whose image in $G$ has order $r$ for some divisor
    $r>1$ of~$d$ .  This gives us the followiung commutative diagram:
    \[
    \begin{diagram}
      \node{L_{4k}(\Q[\Z_d])} \arrow{r}\arrow{s} \node{L_{4k}(\Q\pi)}\arrow{s} \\
      \node{L_{4k}(\Q[\Z_r])} \arrow{r} \node{L_{4k}(\Q G)}
    \end{diagram}
    \]
    Since the class of the form $[1]$ in $L_{4k}(\Q[\Z_r])$ is in the
    image of $L_{4k}(\Q[\Z_d])$, the elements $o\in L_{4k}(\Q G)$
    described above are all in the image of $L_{4k}(\Q\pi)$.
    
    It follows that there are infinitely many $a\in L_{4k}(\pi)$ such
    that $\rho^{(2)}(M,\phi)-\rho^{(2)}(M_a,\phi_a)$ are all nonzero
    and distinct.  From the homology cobordism invariance of
    $\rho^{(2)}(M,\phi)$ (see
    Theorem~\ref{theorem:homology-invariance-der-series}) the desired
    conclusion follows.  
    
\end{proof}

\subsection{Homology cobordism types of 3-manifolds}

In the below statement, $\pi^{(n)}$ denotes the $R$-coefficient
Vogel-Cohn local derived series of a group~$\pi$, where $R=\Z_p$,
$\Z_{(p)}$, or~$\Q$.

\begin{theorem}
  \label{theorem:exotic-homology-cobordism-types-dim3}
  Suppose $M$ is a closed 3-manifold with $\pi=\pi_1(M)$.  If
  $\pi^{(n)}/\pi^{(n+1)}$ is nontrivial for some $n>0$, then there
  exist infinitely many closed 3-manifolds $M_0 = M$, $M_1$, $M_2,
  \ldots$ satisfying the following properties:
  \begin{enumerate}
  \item For any $i$, there is a map $M\to M_i$ which induces
    isomorphisms on $H_*(-;\Z)$ and on $\pi_1(-)/\pi_1(-)^{(k)}$ for
    any~$k$.
  \item For any $i$ and $j$, $M_i$ and $M_j$ have identical Wall
    multisignatures (or equivalently Atiyah-Singer $G$-signatures).
  \item For any $i\ne j$, $M_i$ and $M_j$ are not homology cobordant.
  \end{enumerate}
\end{theorem}

\begin{corollary}
  \label{corollary:homology-cobordism-spherical-space-form}
  For the spherical 3-space form $M=S^3/Q_{8n}$ with fundamental group
  $Q_{8n} = \langle x,y \mid y^{2n} = x^2 = (xy)^2\rangle$ (see
  Milnor~\cite{milnor:1957-2}), there are infinitely many closed
  3-manifolds $M_0=M$, $M_1$, $M_2,\ldots$ such that the $M_i$ are
  homology equivalent to $M$ and have identical Wall multisignatures
  (or equivalently Atiyah-Singer $G$-signatures) and Harvey
  $\L2$-signature invariants $\rho_n$ \cite{Harvey:2006-1}, but no two
  of the $M_i$ are homology cobordant.
\end{corollary}

\begin{proof}
  Due to our previous computation
  (Example~\ref{example:spherical-space-form-groups}), $\pi=Q_{8n}$
  has non-trivial $\pi^{(1)}/\pi^{(2)}$, where $\pi^{(n)}$ denotes the
  $\Z_{(2)}$-coefficient Vogel-Cohn local derived series.  Therefore
  Theorem~\ref{theorem:exotic-homology-cobordism-types-dim3} applies.
  Harvey's invariant $\rho_n$ vanishes for any rational homology
  sphere.
\end{proof}

\begin{remark}
  \leavevmode\Nopagebreak
  \begin{enumerate}
  \item While Harvey's $\L2$-invariants $\rho_n$ cannot distinguish
    the homology cobordism types in
    Corollary~\ref{corollary:homology-cobordism-spherical-space-form},
    our proof illustrates that $\L2$-signatures, combined with our
    injectivity theorem, are effective even for rational homology
    spheres.
  \item Our computation for $\pi=Q_{8n}$ says that
    $\pi^{(n)}/\pi^{(n+1)}$ is a $2$-group for $n\le 1$; it enables an
    alternative proof of
    Corollary~\ref{corollary:homology-cobordism-spherical-space-form}
    using either the result of Levine~\cite{Levine:1994-1}, or the
    $L$-group-valued Hirzebruch-type invariants from iterated
    $2$-covers which were defined and studied by the first
    author~\cite{Cha:2007-1,Cha:2007-2}.
  \end{enumerate}
\end{remark}

\begin{proof}[Proof of
  Theorem~\ref{theorem:exotic-homology-cobordism-types-dim3}]
  Choose a simple closed curve $\eta$ in $M$ representing a nontrivial
  element in $\pi^{(n)}/\pi^{(n+1)}$.  We will use a construction sometimes
  called an {\em infection} along $\eta$ by a knot: identify a tubular
  neighborhood of $\eta$ with $\eta\times D^2$, choose a knot $K$, and
  let $E_K$ be the exterior of $K$.  Let
  \[
  M(\eta,K) = \big(M-(\eta\times \inte D^2) \big) \cup_{\partial} E_K
  \]
  where the pasting map is an orientation reversing homeomorphism on
  the boundary such that $\eta\times \{*\},\ \{*\}\times S^1 \subset
  \eta\times S^1$ are identified with a meridian and longitude of $K$
  on~$\partial E_K$.  We will show that the 3-manifolds $M_i$ can be
  obtained in this way by choosing the knots $K$ appropriately.

  Let $K_0$ be the unknot.  Then there is a homology equivalence $E_K
  \to E_{K_0}$ that extends a homeomorphism on the boundary preserving
  the periphral structure.  Pasting this with the identity map on
  $M-(\eta\times \inte D^2)$, we obtain a homology equivalence
  $M(\eta,K) \to M(\eta,K_0)=M$.

  Let $G=\pi_1(M(\eta,K))$.  Then the induced map $G\to \pi$ is
  2-connected on homology (for any~$R$).  Therefore by the
  injectivity theorem~\ref{theorem:injectivity}, the induced map $G/G^{(k)} \to \pi/\pi^{(k)}$
  is injective for any~$k$.  Since $\pi_1(E_K) \to \pi_1(E_{K_0})$ is
  surjective, $G\to \pi$ is surjective.  It follows that $G/G^{(k)}
  \to \pi/\pi^{(k)}$ is an isomorphism for any~$k$.

  Since $n\ge 1$, $[\eta]=0$ in $H_1(M)$.  From this it follows that
  Wall's multisignatures of $M(\eta,K)$ are identical with those of~$M$
  (e.g., see \cite{Friedl:2003-1} and Section 5 of \cite{Cha:2007-1}).

  In order to distinguish the homology cobordism classes of the
  $M(\eta,K)$, we will use $\L2$-signatures.  Consider
  $\rho^{(2)}(M(\eta,K),\psi)$, where $\psi$ is the quotient map $G
  \to G/G^{(n+1)} \cong \pi/\pi^{(n+1)}$.  
  
By previously known arguments (e.g., see~\cite{Cochran-Orr-Teichner:2002-1},~\cite{Cochran-Harvey-Leidy:2009}),
we have
  \[
  \rho^{(2)}(M(\eta,K),\psi) = \rho^{(2)}(M,\phi) +
  \rho^{(2)}(M_K,\alpha)
  \]
  where $\phi$ is the quotient map $\pi \to \pi/\pi^{(n+1)}$, $M_K$ is
  the zero-surgery manifold of $K$, and $\alpha\colon \pi_1(M_K) \to
  \pi/\pi^{(n+1)}$ is the map induced by the restriction of $\psi$
  on~$E_K$.

  Let $d$ be the order of $[\eta]$ in $\pi/\pi^{(n+1)}$.  ($d$ may be
  $\infty$.)  Since $\pi^{(n)}$ is a normal subgroup and $\pi_1(M_K)$
  is normally generated by a meridian of $K$ which is identified with
  a parallel copy of $\eta$, $\Im \alpha$ is contained in
  $\pi^{(n)}/\pi^{(n+1)}$ which is an abelian group.  It follows that
  $\Im \alpha$ is isomorphic to $\Z_d$ (where $\Z_\infty$ is
  understood as the infinite cyclic group if $d=\infty$).  Therefore
  by the induction property, $\rho^{(2)}(M_K,\alpha) =
  \rho^{(2)}(M_K,\beta)$ where $\beta$ is the surjection $\pi_1(M_K)
  \to \Im \alpha \cong \Z_d$.  Appealing to
  Lemma~\ref{lemma:$\L2$-signature-of-knots} below, we can choose
  infinitely many knots $K_0$ $(=$ unknot$)$, $K_1$, $K_2,\ldots$ such
  that the values $\rho^{(2)}(M_{K_i},\beta)$ are mutually distinct.
  It follows that the $\rho^{(2)}(M(K_i,\eta),\phi_{K_i})$ are all
  distinct.  By the homology cobordism invariance of
  $\rho^{(2)}(M(K_i,\eta),\psi)$
  (Theorem~\ref{theorem:homology-invariance-der-series}), the
  3-manifolds $M_i = M(\eta, K_i)$ are not homology cobordant.
\end{proof}

\begin{lemma}
  \label{lemma:$\L2$-signature-of-knots}
  \leavevmode\Nopagebreak
  Suppose $K$ is a knot in $S^3$ such that zero-surgery on $K$ yields the three manifold~$M$.  Let $A$ be a Seifert matrix for~$K$.  Let
  \[
  \sigma_K(\omega)=\sign\big((1-\omega)A + (1-\omega^{-1})A^T\big)
  \]
  be the Levine-Tristram signature function of $K$ which is defined
  for $\omega\in S^1 \subset \C$.
  \begin{enumerate}
  \item \cite[Proposition~5.1]{Cochran-Orr-Teichner:2002-1} For the
    abelianization map $\phi\colon \pi_1(M) \to \Z$, we have
    \[
    \rho^{(2)}(M,\phi)=\int_{S^1} \sigma_K(\omega)\, d\omega
    \]
    where the integral is over $S^1$ normalized to unit length.
  \item Suppose $\phi_d\colon \pi_1(M) \to \Z_d$ is a surjection.
    Then
    \[
    \rho^{(2)}(M,\phi_d)=\frac{1}{d}\cdot \sum_{k=0}^{d-1}
    \sigma_K(\zeta_d^k)
    \]
    where $\zeta_d=e^{2\pi\sqrt{-1}/d}$ is the $d$th primitive root of
    unity.
  \end{enumerate}
\end{lemma}

\begin{proof}
 Construct in the usual way, a compact 4-manifold $W$ with
  boundary $M$ over $\Z$ whose $\Z[t,t^{-1}]$-coefficient intersection
  form is $(1-t)A+(1-t^{-1})A^T$.  (1)~is proved in
  \cite{Cochran-Orr-Teichner:2002-1} by computing the $\L2$-signature
  of the intersection form of~$W$.  One can prove (2) similarly, as
  follows: let $G=\Z_d$.  Since $G$ is finite, the von
  Neumann group ring $\N G$ is equal to the ordinary group ring
  $\C G$, and $d\cdot \ldim M = \dim_{\C} M$ for any
  $\N G$-module~$M$.  Viewing the bordism $W$ described above as a
  bordism over $G$, the $\N G$-coefficient intersection form
  is exactly
  \[
  (1-g)A + (1-g^{-1})A^T
  \]
  where $g$ is the generator of $G$ which is the image of the
  positive meridian.  It can also be seen that the ordinary signature
  of $W$ is equal to $\sigma_K(1)$, which is always zero.  From these
  observations (2) follows.
\end{proof}

\appendix

\section{Remarks on module and ring localizations}

\subsection{Two definitions of the Bousfield module localization}

In \cite{Bousfield:1975-1}, Bousfield defined and studied localization
of $Z\pi$-modules with respect to the class $\HZ$ of $Z\pi$-module
homomorphisms $\alpha\colon A\to B$ such that the induced map
$H_i(\pi;A) \to H_i(\pi;B)$ is an isomorphism for $i=0$ and a
surjection for $i=1$.  His arguments readily extend to the case of
a group $\pi$ over another group~$\Gamma$.  Namely, there is a
localization functor with respect to the class, which will be denoted
by $\HZ$ as well, of $R\pi$-module homomorphisms $\alpha\colon A\to B$
that induces an isomorphism $A\otimesover{R\pi} R\Gamma \to
B\otimesover{R\pi} R\Gamma$ and a surjection $\Tor_1^{R\pi}(R\Gamma,A)
\to \Tor_1^{R\pi}(R\Gamma,B)$.

Our goal is to show that the Bousfield localization with respct to
$\HZ$ is equal to the localization with respect to the following
class: let $\W$ be the collection of $R\pi$-module homomorphisms
$\alpha\colon F \to F'$ such that $F$ and $F'$ are $R\pi$-free and
$\alpha\otimes 1\colon F\otimesover{R\pi} R\Gamma \to
F'\otimesover{R\pi} R\Gamma$ is an isomorphism.  Although surely well-known, we could
find no proof in the literature for the following, and provide the proof here.

\begin{theorem}
  A $R\pi$-module $M$ is local with respect to $\HZ$ if and only if
  $M$ is local with respect to~$\W$.
\end{theorem}  

The following is an immediate consequence:

\begin{corollary}
  The localization with respect to $\HZ$ is equal to the localization
  with respect to~$\W$.
\end{corollary}

\begin{proof}[Proof of Theorem]
  Observe that the class $\W$ is contained in~$\HZ$.  It follows that
  if $M$ is local with respect to $\HZ$, then $M$ is local with
  respect to~$\W$.

  To prove the converse, suppose $M$ is local with respect to $\W$,
  and consider a diagram
  \[
  \begin{diagram}
    \node{A} \arrow{e,t}{\alpha}\arrow{s}
    \node{B}
    \\
    \node{M}
  \end{diagram}
  \]
  with $\alpha$ in $\HZ$.  We have to show that there is a unique map
  $B\to M$ making the diagram commute.

  Let $N=\Ker\{\pi \to \Gamma\}$ and $I=\Ker\{R\pi \to R\Gamma\}$.
  Choose a set of generators $\{b_i\}_i$ of $B$.  Since
  \[
  A\otimesover{R\pi}R\Gamma \cong B\otimesover{R\pi}R\Gamma \cong
  B/\langle gb-b \mid g\in N, b\in B\rangle
  \]
  we can write each $b_i$ as
  \[
  b_i = \alpha(a_i) + \sum_j r_{ij} b_j
  \]
  where $a_i\in A$ and $r_{ij} \in I$.  We consider the following
  system of equations with variable $x_i$ (one variable for each
  generator $b_i$):
  \[
  S=\Big\{x_i = a_i + \sum_j r_{ij} x_j \Big\}_i
  \]
  When $S$ is a system of equation in $A$ as above, we define $A_S$ to
  be the module obtained from $A$ by adjoining the variables $x_i$ as
  additional generators satisfying the relations given by the
  equations in $S$.  Precisely,
  \[
  A_S = (A \oplus F\langle x_i \rangle) \Big/ \Big\langle x_i-a_i-\sum_j
  r_{ij}x_j \Big\rangle
  \]
  where $F\langle x_i \rangle$ denotes the free $R\pi$-module
  generated by the~$x_i$.  Note that the natural map $A\to A_S$
  induces an isomorphism $A\otimesover{R\pi} R\Gamma \to
  A_S\otimesover{R\pi} R\Gamma$ provided $r_{ij}\in I$.  Also, in our
  case, there is a well-defined map $A_S \to B$ sending $x_i$
  to~$b_i$.

  Choose a surjection $F\to A$ of a free $R\pi$-module~$F$.  Let $S'$
  be a system of equations in $F$ which is a lift of $S$, that is,
  $S'$ is obtained from $S$ by replacing each $a_i$ by a pre-image of
  $a_i$ in~$F$.  We have the following commutative diagram:
  \[
  \begin{diagram}
    \node{F} \arrow{e} \arrow{s}
    \node{F_{S'}} \arrow{s} \arrow{ssw,..}
    \\
    \node{A} \arrow{e} \arrow{s}
    \node{A_S} \arrow{e}
    \node{B}
    \\
    \node{M}
  \end{diagram}
  \]
  Since $F\to F_{S'}$ is in $\W$, there is $F_{S'} \to M$ making the
  diagram commute.

  It can be checked that $\Ker\{F_{S'} \to A_S\}$ is equal to the
  image of $\Ker{F\to A}$ under $F\to F_{S'}$.  Therefore, $F_{S'}\to
  M$ induces $A_S \to M$.  Observe the following facts:
  \begin{enumerate}
  \item $A_S \otimesover{R\pi} R\Gamma \cong A \otimesover{R\pi}
    R\Gamma \cong B\otimesover{R\pi} R\Gamma$.
  \item $\Tor_1^{R\pi}(A_S,R\Gamma) \to \Tor_1^{R\pi}(B,R\Gamma)$ is
    surjective, since so is $\Tor_1^{R\pi}(A,R\Gamma) \to
    \Tor_1^{R\pi}(B,R\Gamma)$.
  \end{enumerate}
  Let $K=\Ker\{A_S \to B\}$.  Looking at the $\Tor$ long exact
  sequence obtained from $0\to K \to A_S \to B\to 0$, from the above
  observations it follows that $R\Gamma\otimesover{R\pi}K=0$.  By the
  lemma stated and proved below, it follows that the image of $K$ in
  $M$ under $A_S\to M$ is zero.  Therefore, from (1), it follows that
  there is an induced map $B\to M$.  The map $B\to M$ is uniquely
  determined since $F_{S'} \to M$ is unique.
\end{proof}

\begin{definition}
  For $\pi$ over $\Gamma$, an $R\pi$-module $N$ is \emph{perfect} if
  $N\otimesover{R\pi}R\Gamma=0$.
\end{definition}

\begin{lemma}
  If $M$ is local with respect to $\W$, then $M$ has no nontrivial
  perfect submodules.  Consequently any map of a perfect module into
  $M$ is zero.
\end{lemma}

\begin{proof}
  Suppose $N$ is a perfect submodule of~$M$.  Choose generators
  $\{a_i\}$ of~$N$.  Since $N \otimesover{R\pi} R\Gamma = N/\langle
  ga-a \mid g\in \Ker\{\pi\to\Gamma\}, a\in N \rangle = 0$, we can
  write each $a_i$ as $ a_i = \sum_j r_{ij}a_j $ where $r_{ij}\in
  I=\Ker\{R\pi \to R\Gamma\}$.  Let $F\langle x_i\rangle$ be the free
  $R\pi$-module generated by variables $x_i$ (one variable $x_i$ for
  each $b_i$) as before, and define $\alpha\colon F\langle x_i\rangle
  \to F\langle x_i\rangle$ by $\alpha(x_i) = x_i-\sum_j r_{ij}x_j$.
  Consider
  \[
  \begin{diagram}
    \node{F\langle x_i \rangle} \arrow{e,t}{\alpha} \arrow{s,l}{0}
    \node{F\langle x_i \rangle} \arrow{sw,b,..}{\phi}
    \\
    \node{M}
  \end{diagram}
  \]
  Observe that $\alpha$ is in~$\W$.  The map $\phi\colon F\langle x_i
  \rangle \to M$ given by $\phi(x_i)=a_i$ makes the diagram commute.
  Instead of $\phi$, the zero map also makes it commute.  Therefore,
  since $M$ is local with respect to $\W$, $\phi=0$ by the uniqueness.
  It follows that $a_i=0$.
\end{proof}

\subsection{Ring localization as module localization}

Suppose $\pi$ is a group, $R$ is a commutative ring with unity,
$\Omega$ is a collection of $R\pi$-module homomorphisms, and $E$ is a
localization functor on the category of $R\pi$-modules with respect
to~$\Omega$, which is endowed with a natural transformation $M \to
E(M)$.  Our goal is to give a proof that $E(\Z\pi)$ is a ring
localization in the sense of
Section~\ref{subsection:ring-localization}.

\begin{theorem}
  There is a ring structure on $\Lambda=E(R\pi)$ satisfying the following:
  \begin{enumerate}
  \item The map $R\pi \to \Lambda$ is a ring homomorphism.
  \item For any $R\pi$-module $M$, there is a natural isomorphism
    $E(M) \cong M\otimesover{R\pi} \Lambda$.
  \item For any $\alpha\colon A \to B$ in $\Omega$, $\alpha\otimes
    1\colon A\otimesover{R\pi} \Lambda \to B\otimesover{R\pi}
    \Lambda$ is an isomorphism.
  \item $R\pi \to \Lambda$ is initial among objects satisfying (3),
    that is, if a ring homomorphism $R\pi\to \Lambda'$ satisfies (3),
    then there is a unique ring homomorphism $\Lambda \to \Lambda'$
    making the following diagram commute:
    \[
    \begin{diagram}
      \node{R\pi} \arrow{e}\arrow{s}
      \node{\Lambda} \arrow{sw,..}
      \\
      \node{\Lambda'}
    \end{diagram}
    \]
  \end{enumerate}
\end{theorem}

\begin{proof}
  (1) Denote by $i$ the map $R\pi \to E(R\pi)=\Lambda$.  The ring
  structure on $\Lambda$ is defined as follows.  Since $R\pi
  \otimesover{R\pi} \Lambda \cong \Lambda$ and the functor $E$
  commutes with direct sum, for any free $R\pi$-module $F$, we have
  the following commutative diagram:
  \[
  \begin{diagram} 
    \node{\textstyle F\otimesover{R\pi} R\pi}
    \arrow{e,t}{\cong}\arrow{s,l}{1\otimes i}
    \node{F} \arrow{s}
    \\
    \node{\textstyle F\otimesover{R\pi}\Lambda} \arrow{e,t}{\cong}
    \node{E(F)}
  \end{diagram}
  \]
  Choose a free resolution $F_0 \to F_1 \to \Lambda$ over~$R\pi$.
  Then we have
  \[
  \begin{diagram}
    \node{\textstyle F_0 \otimesover{R\pi} \Lambda} \arrow{e} \arrow{s}
    \node{\textstyle F_0 \otimesover{R\pi} \Lambda} \arrow{e} \arrow{s}
    \node{\textstyle \Lambda \otimesover{R\pi} \Lambda} \arrow{e}
    \node{0}
    \\
    \node{E(F_0)} \arrow{e}
    \node{E(F_1)} \arrow{e}
    \node{E(\Lambda)} \arrow{e}
    \node{0}
  \end{diagram}
  \]
  Rows are exact since tensoring and $E$ are right exact, and vertical
  arrows are isomorphisms.  Since $E$ is an idempotent, it follows
  that there is an isomorphism
  \[
  m\colon \Lambda\otimesover{R\pi} \Lambda \to \Lambda= E(\Lambda)
  \]
  which gives ring multiplication on $\Lambda$.  Obviously $R\pi \to
  \Lambda$ is a ring homomorphism.

  (2) For a given $R\pi$-module $M$, choose a presentation $F_0 \to
  F_1 \to M \to 0$ with $F_0$, $F_1$ free.  Replacing the resolution
  of $\Lambda$ in the above argument of (1) by the resolution of $M$,
  we obtain an isomorphism $M\otimesover{R\pi}\Lambda \to E(M)$.  It
  can be verified that this isomorphism is independent of the choice
  of the presentation.

  (3) Any $\alpha\colon A \to B$ in $\Omega$ induces an isomorphism
  $E(A) \to E(B)$ since $E$ is a localization with respect
  to~$\Omega$.  Since $E(A) \cong A\otimesover{R\pi}\Lambda$, the
  desired conclusion follows.

  (4) Suppose a ring homomorphism $j\colon R\pi \to \Lambda'$
  satisfies~(3).  First we will show that $\Lambda'$ is a local
  $R\pi$-module with respect to~$\Omega$.  Suppose $\alpha\colon A \to
  B$ in $\Omega$ and $\phi\colon A \to \Lambda'$ is given.
  
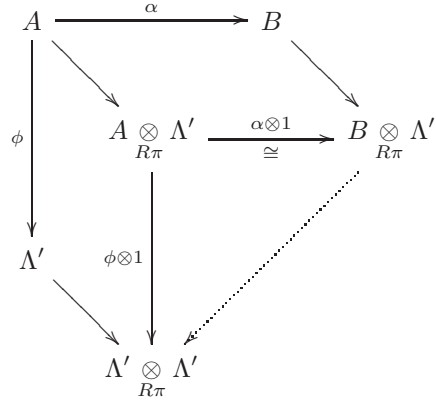
\[
  \begin{diagram}
    \node{A} \arrow[2]{e,t}{\alpha} \arrow{se} \arrow[2]{s,l}{\phi}
    \node[2]{B} \arrow{se}
    \\
    \node[2]{A\otimesover{R\pi} \Lambda'} \arrow[2]{s,l}{\phi\otimes 1}
    \arrow[2]{e,tb}{\alpha\otimes 1}{\cong}
    \node[2]{B\otimesover{R\pi} \Lambda'} \arrow[2]{sw,..}
    \\
    \node{\Lambda'} \arrow{se}
    \\
    \node[2]{\Lambda' \otimesover{R\pi} \Lambda'}
  \end{diagram}
  \]
  In the above diagram, $\alpha\otimes 1$ is an isomorphism by~(2).
  Let $m'\colon \Lambda'\otimesover{R\pi}\Lambda' \to \Lambda'$ be the
  ring multiplication and let $\psi\colon B \to \Lambda'$ be the
  composition of $B\to B\otimesover{R\pi} \Lambda'$ and $m'(\phi\times
  1)(\alpha\otimes 1)^{-1}$.  Then $\psi\alpha=\phi$.  If
  $\psi'\alpha=\phi$ for some $\psi'\colon B \to \Lambda'$, then
  $\psi'\otimes 1\colon B\otimesover{R\pi}\Lambda' \to
  \Lambda'\otimesover{R\pi}\Lambda'$ should be equal to $(\phi\otimes
  1)(\alpha\otimes 1)^{-1}$, being the unique map making the triangle
  commute.  From this it follows that $\phi'=\phi$ by a diagram chase.
  This proves that the $R\pi$-module $\Lambda'$ is local.

  By the universal property of the module localization
  $\Lambda=E(R\pi)$, there is a unique $R\pi$-module homomorphism
  $f\colon \Lambda \to \Lambda'$ making the following diagram commute:
  
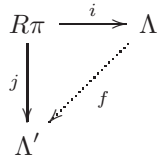
\[
  \begin{diagram}
    \node{R\pi} \arrow{e,t}{i}  \arrow{s,l}{j}
    \node{\Lambda}\arrow{sw,b,..}{f}
    \\
    \node{\Lambda'}
  \end{diagram}
  \]
  We claim that $f$ is a ring homomorphism.  To prove this, consider
  the following diagram, where $m$, $m_0$, $m'$ are ring
  multiplication.
  
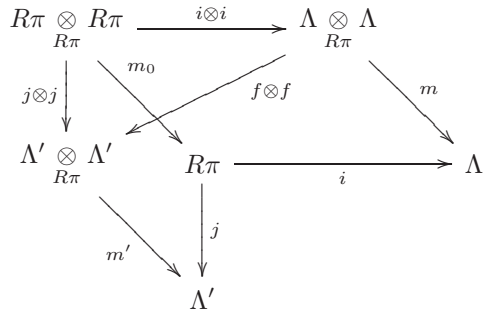
\[
  \begin{diagram}
    \node{R\pi\otimesover{R\pi} R\pi}
    \arrow[2]{e,t}{i\otimes i}
    \arrow{s,l}{j\otimes j}
    \arrow{se,t,1}{m_0}
    \node[2]{\Lambda\otimesover{R\pi}\Lambda} 
    \arrow{se,t}{m} \arrow{sww,b,1}{f\otimes f}
    \\
    \node{\Lambda'\otimesover{R\pi}\Lambda'} \arrow{se,b}{m'}
    \node{R\pi} \arrow[2]{e,b}{i} \arrow{s,r}{j}
    \node[2]{\Lambda}
    \\
    \node[2]{\Lambda'}
  \end{diagram}
  \]
  It can be seen that the diagram commutes.  Now look at $m'(f\otimes
  f) m^{-1}\colon \Lambda \to \Lambda'$.  Then by a diagram chase, one
  can verify 
  \[
  m'(f\otimes f) m^{-1}i = m'(f\otimes f) (i\otimes i) m_0^{-1} =
  m'(j\otimes j) m_0^{-1} = j.
  \]
  From the uniqueness of $f$, it follows that $m'(f\otimes
  f)m^{-1}=f$, that is, $f$ is a ring homomorphism.
\end{proof}

\bibliographystyle{amsplainabbrv}
\renewcommand{\MR}[1]{}

\bibliography{research}

\end{document}